\documentclass[11pt] {article} 
\usepackage[backref]{hyperref}
  \usepackage{amsmath}
    \usepackage{amssymb}

\newtheorem{theorem}{Theorem}[section]
\newtheorem{lemma}{Lemma}[section]

\newcommand{\eqnsection}{
   \renewcommand{\theequation}{\thesection.\arabic{equation}}
   \makeatletter
   \csname @addtoreset\endcsname{equation}{section} 
   \makeatother}
 
\def \be{\begin{equation}}
\def \ee{\end{equation}}
\def \bt{\begin{theorem}} 
\def \et{\end{theorem}} 
\def \bl{\begin{lemma}} 
\def \el{\end{lemma}}
\def \bea{\begin{eqnarray}}
\def \eea{\end{eqnarray}}
\def \bas{\begin{eqnarray*}}
\def \eas{\end{eqnarray*}}

\def \al{\alpha}  
\def \bb{\beta}
\def \ga{\gamma} 
 
\def \de{\delta} 
\def \De{\Delta} 
\def \ep{\epsilon}

\def \la{\lambda}  
\def \La{\Lambda}

\def \om{\omega}
\def \Om{\Omega}

\def \si{\sigma}

\def \th{\theta}

\def \ze{\zeta}

\def \ff{\infty}
\def \wh{\widehat}
\def \wt{\widetilde}

\def \rar{\rightarrow}

\def \PP{{\cal P}}

\def \({\left(}
\def \){\right)}

\def \lc{\left\{}
\def \rc{\right\}}

\def \nn{\nonumber}

\def \bc{\begin{center} }
\def \ec{\end{center} }
\def \bs{\begin{slide} }
\def \es{\end{slide} }

\def\square{{\vcenter{\vbox{\hrule height.3pt
        \hbox{\vrule width.3pt height5pt \kern5pt
           \vrule width.3pt}
        \hrule height.3pt}}}}
\def\qed{{\hfill $\square$ \bigskip}}

 \def \Rev({\mbox{\rm Rev}(}

\eqnsection
\begin{document}

\title{Lectures on Isomorphism Theorems}

 \author{Jay Rosen \thanks{Research   was supported by  grants from the National Science Foundation.  }}
\maketitle

 
 \newpage
 \tableofcontents

\bibliographystyle{amsplain}

\section{Introduction}

These notes originated in a series of lectures I gave in Marseille in May, 2013. I was invited to give an introduction to the isomorphism theorems, originating with Dynkin, \cite{Dynkin83},  \cite{Dynkin84}, which connect Markov local times and Gaussian processes. This is an area I had worked on some time ago, and even written a book about, \cite{book}, but had then moved on to other things. However,  isomorphism theorems have  become of interest once again, both because of   new applications to the study of  cover times of graphs and Gaussian fields, \cite{DLP, D1, D2, DZ1, DZ2}, and because of new isomorphism theorems for non-symmetric Markov processes and their connection with loop soups and Poisson processes, \cite{EK, FR, Le Jan, Le Jan1, LMR, LMR2}. Thus I felt the time was ripe for a new introduction to this topic.

I greatly enjoyed giving these lectures, since I felt free to focus on what I consider to be the basic ideas. Writing my book with Marcus took a lot of time and effort since we wanted to make sure that all the details were carefully explained. In these notes I have tried to preserve the informal atmosphere of the lectures, and often simply refer the reader to the book \cite{book} and other sources for details. 

The actual lectures covered the material which appears in sections 2-7. This begins with some introductory material on Gaussian processes and Markov processes and then studies in turn the isomorphism theorems of Dynkin, Eisenbaum and the generalized second Ray-Knight theorem. In each case we give a proof and a sample application. We then introduce loop soups and permanental processes, the ingredients we use to develop an  isomorphism theorem for non-symmetric Markov processes. Along the way we gain new insight into the reason that Gaussian processes appear in the isomorphism theorems in the symmetric case. Chapters 8-10 contain the material I would have liked to include in the lectures, but had to skip because of lack of time. Having developed some general material on Poisson processes in Section 7, we make use of it in the next two sections. Section 8 contains an excursion theory proof of the generalized second Ray-Knight theorem, and section 9 explains a similar theorem for random interlacements. Up till this point, the proofs I give use the method of moments, which for me is the simplest and clearest way to prove isomorphism theorems. In section 10 we explain how to prove these theorems using the method of Laplace transforms.  Some may prefer this approach because it is more `automatic' and doesn't involve the sometimes subtle combinatorics which come up when dealing with moments.

\section{Gaussian processes}\label{sec-gp}

A real valued random variable $X$ is called a Gaussian random variable if  
\begin{equation}
E\(e^{i\la X}\)=e^{im\la -\si^{2}\la^{2}/2}, \hspace{.2 in}\forall \la\in R^{1}\label{gp.1}
\end{equation}
for some numbers $m,\si$. Differentiating in $\la$ we see that
\begin{equation}
E(X)=m, \hspace{.2 in}V(X)=\si^{2}.\label{gp.2}
\end{equation}
We can always eliminate the $m$ by subtracting it from $X$. From now on we assume that 
$E(X)=0$, so that (\ref{gp.1}) becomes
\begin{equation}
E\(e^{i\la X}\)=e^{-E((\la X)^{2})/2}, \hspace{.2 in}\forall \la\in R^{1}.\label{gp.3}
\end{equation}

A  random vector $X=(X_{1}, \ldots, X_{n})\in R^{n}$ is called a Gaussian random vector if  
for each $y\in R^{n}$, $(y,X)$ is a Gaussian random variable. Thus we have
\begin{equation}
E\(e^{i(y,X)}\)=e^{-E((y,X)^{2})/2}, \hspace{.2 in}\forall y\in R^{n}.\label{gp.4}
\end{equation}
The $n\times n$ matrix $C=\{C_{i,j},1\leq i,j\leq n \}$ with entries 
\begin{equation}
C_{i,j}=E(X_{i}X_{j}), \hspace{.2 in}1\leq i,j\leq n\label{gp.5}
\end{equation}
is called the covariance matrix of $X$, and we can write
\begin{equation}
(y,X)^{2}=\sum_{i,j=1}^{n}C_{i,j}y_{i}y_{j}=(y,Cy),\label{gp.6}
\end{equation}
so that (\ref{gp.4}) can be written as
\begin{equation}
E\(e^{i(y,X)}\)=e^{-(y,Cy)/2}, \hspace{.2 in}\forall y\in R^{n}.\label{gp.7}
\end{equation}
It follows from (\ref{gp.5}) that $C$ is symmetric and from (\ref{gp.6}) that $C$ is positive definite.
We now show that conversely, any symmetric positive definite  $n\times n$   matrix $B$ is the covariance matrix of some Gaussian random vector in $R^{n}$.

To see this, we first note that if $A$ is any $p\times n$ matrix and we write $Z=AX$, then by (\ref{gp.7}) we have 
\begin{equation}
E\(e^{i(y,Z)}\)=E\(e^{i(A^{t}y,X)}\)=e^{-(A^{t}y,CA^{t}y)/2}=e^{-(y,ACA^{t}y)/2}, \hspace{.2 in}\forall y\in R^{p}.\label{gp.8}
\end{equation}
so that $Z=AX$ is Gaussian random vector with covariance matrix $ACA^{t}$.

If  $B$ is a symmetric positive definite  $n\times n$ matrix, then there exists    a symmetric  matrix $A$ with $B=A^{2}$. To see this recall that any symmetric matrix is diagonalizable, so we can find an orthonormal system of vectors $u_{i}, 1\leq i\leq n$ such that $Bu_{i}=\la_{i}u_{i}, 1\leq i\leq n$, and the fact that $B$ is   positive definite implies that all $\la_{i}\geq 0$. We can then define the matrix $A$ by setting $Au_{i}=\la^{1/2}_{i}u_{i}, 1\leq i\leq n$. If we now take $X$ to be a vector whose components are independent standard normals, so that the covariance matrix of $X$ is $I$,
it follows from the above that $Z=AX$ is a Gaussian random vector with covariance matrix $B$.

If   $S$ is a general set, a stochastic process $G=\{G_{x}, x\in S\}$ is called a Gaussian process on $S$ if for any $n$ and any  $x_{1}, \ldots,x_{n}\in S$, $(G_{x_{1}}, \ldots, G_{x_{n}})$ is a Gaussian random vector. Then the function
\begin{equation}
C(x,y)=E(G_{x}G_{y}), \hspace{.2 in}x,y\in S\label{gp.9}
\end{equation} 
on $S\times S$ is called the covariance function of $G$. Using the above and Kolmogorov's extension theorem we see that there is a correspondence between Gaussian processes on $S$
and symmetric  positive definite functions on $S\times S$.

Example: Let $S=R_{+}^{1}$, and let $C(s,t)=s\wedge t=\int 1_{[0,s]}(x)1_{[0,t]}(x)\,dx$. Then $C(s,t)$ is positive definite since 
\begin{equation}
\sum_{i,j=1}^{n}C(s_{i},s_{j})y_{i}y_{j}=\int \(\sum_{i=1}^{n} y_{i}1_{[0,s_{i}]}(x)\)^{2}\,dx\geq 0,\label{gp.10}
\end{equation}
so there exists a   Gaussian process $B=\{B_{s}, s\in R_{+}^{1}\}$. Note that is $s<t<t'$ we have 
\begin{equation}
E\(B_{s}(B_{t'}-B_{t})\)=s\wedge t-s\wedge t'=0\label{gp.11}
\end{equation}
so that $B$ has orthogonal increments, which are then independent by (\ref{gp.4}). Hence $B$ is `almost' Brownian motion. What is missing is a continuous version, which can be established in the usual ways.

\subsection{Gaussian moment formulas}\label{sec-gmf}

Let $G=\{G_{x}, x\in S\}$ be a Gaussian process with covariance function $C$. We present several 
Gaussian moment formulas which will be used to prove our Isomorphism Theorems. The basic formula is 

\begin{equation}
E\(\prod_{i=1}^{n}G_{x_{i}}\)=\sum_{p\in\mathcal{R}_{n}}\prod_{(i_{1},i_{2})\in p}C(x_{i_{1}},  x_{i_{2}})\label{gp.12}
\end{equation}
where $\mathcal{R}_{n}$ denotes the set of pairings $p$ of the indices $[1,n]$, and the product runs over all pairs in $p$. In particular, this is empty when $n$ is odd, in which case the left hand side is zero by symmetry.

Proof: We write (\ref{gp.7}) as 
\begin{equation}
E\(e^{i\sum_{j=1}^{n}z_{j} G_{x_{j}}    }\)=e^{-\sum_{j,k=1}^{n}z_{j}z_{k}C(x_{j},  x_{k})/2}.\label{gp.13}
\end{equation}
We differentiate successively in $z_{1},\ldots,z_{n}$ and after differentiating in $z_{j}$ we set 
$z_{j}=0$.

To begin, we differentiate in $z_{1}$ and then  set 
$z_{1}=0$, to obtain
\begin{equation}
iE\(G_{x_{1}}e^{i\sum_{j=2}^{n}z_{j} G_{x_{j}}    }\)=\(-\sum_{k=2}^{n}z_{k}C(x_{1},  x_{k})\)e^{-\sum_{j=2}^{n}z_{j} C(x_{j},  x_{k})/2}.\label{gp.13}
\end{equation}
We then differentiate in $z_{2}$, using the product rule for the right hand side, and after    setting 
$z_{2}=0$ we obtain
\begin{eqnarray}
&&-E\(G_{x_{1}}G_{x_{2}}e^{i\sum_{j=3}^{n}z_{j} G_{x_{j}}    }\)=-C(x_{1},  x_{2}) e^{-\sum_{j=3}^{n}z_{j} C(x_{j},  x_{k})/2}
\label{gp.14}\\
&&\hspace{.4 in}+ \(-\sum_{k=3}^{n}z_{k}C(x_{1},  x_{k})\)\(-\sum_{k=3}^{n}z_{k}C(x_{2},  x_{k})\)e^{-\sum_{j=3}^{n}z_{j} C(x_{j},  x_{k})/2}.  \nonumber
\end{eqnarray}
By now it should be clear that by continuing this process we obtain (\ref{gp.12}).\qed

Our next formula is:
\begin{equation}
E\(\prod_{i=1}^{n}G^{2}_{x_{i}}\)=\sum_{A_{1}\cup \cdots\cup A_{j}=[1,n]}\prod^{j}_{l=1}2^{| A_{l}|-1}\mbox{cy}(A_{l}),\label{gp.15}
\end{equation}
where the sum is over all (unordered)  partitions $A_{1}\cup \cdots\cup A_{j}$ of $  [1,n]$ and, if we have $A_{l}=\{l_{1},l_{2},\cdots, l_{| A_{l}|}\}$ then the cycle function $\mbox{cy}(A_{l})$ is defined as 
\begin{equation}
\mbox{cy}(A_{l})=\sum_{\pi\in \mathcal{P}^{\odot}_{| A_{l}|}}C(x_{l_{\pi(1)}},  x_{l_{\pi(2)}})\cdots C(x_{l_{\pi (| A_{l}|)}},  x_{l_{\pi(1)}}),\label{gp.16}
\end{equation}
where $\mathcal{P}_{k}^{\odot}$ denotes the set of permutations of $[1,k]$ on the circle.   (For example,  $(1,2,3)$,  $(3,1,2)$ and $(2,3,1)$ are considered to be the same permutation   $\pi\in \mathcal{P}^{\odot}_{3 }$.)

Proof: On the left hand side of (\ref{gp.15}) each $G_{x_{i}}$ appears twice. We can arbitrarily consider one of the two $G_{x_{i}}$'s as the `red' $G_{x_{i}}$ and the other as the `green' $G_{x_{i}}$. Consider first the `red' $G_{x_{1}}$. By (\ref{gp.12}) it is paired with some  $G_{x_{i}}$. If it is paired with the `green' $G_{x_{1}}$, we set $A_{1}=\{1\}$, in which case $\mbox{cy}(A_{1})=C(x_{1},x_{1})$. Otherwise, the `red' $G_{x_{1}}$ is coupled with one of the two $G_{x_{j}}$'s for some $j\neq 1$, giving a factor of $2C(x_{1},x_{j})$ and so we continue until eventually we are paired with the `green' $G_{x_{1}}$ which gives the factor $C(x_{\cdot},x_{1})$. $A_{1}$ consists of those  $j$ such that $G_{x_{j}}$ has been used.  Beginning again  with some  $G_{x_{j}}$ not used yet and iterating  we are led to (\ref{gp.15}).\qed

For later reference it will be useful to write (\ref{gp.15}) as
\begin{equation}
E\(\prod_{i=1}^{n}G^{2}_{x_{i}}/2\)=\sum_{A_{1}\cup \cdots\cup A_{j}=[1,n]}\prod^{j}_{l=1}{1 \over 2}\,\mbox{cy}(A_{l}).\label{gp.17}
\end{equation}

Our last formula for now is:
\begin{equation}
E\(G_{a}G_{b}\prod_{i=1}^{n}G^{2}_{x_{i}}/2\)=\sum_{A\subseteq [1,n]}\mbox{ch}(A; a,b)\sum_{A_{1}\cup \cdots\cup A_{j}=[1,n]-A}\prod^{j}_{l=1}{1 \over 2}\,\mbox{cy}(A_{l}),\label{gp.18}
\end{equation}
where the sum is over all (unordered)  partitions $A_{1}\cup \cdots\cup A_{j}$ of $  [1,n]-A$ and, if   $A=\{l_{1},l_{2},\cdots, l_{| A|}\}$ then the chain function $\mbox{ch}(A; a,b)$ is defined as 
\begin{equation}
\mbox{ch}(A; a,b)=\sum_{\pi\in \mathcal{P}_{| A|}}C(x_{a},  x_{l_{\pi(1)}})C(x_{l_{\pi(1)}},  x_{l_{\pi(2)}})\cdots C(x_{l_{\pi (| A_{l}|)}},  x_{b}),\label{gp.19}
\end{equation}
where $\mathcal{P}_{k}$ denotes the set of permutations of $[1,k]$.
Using (\ref{gp.17}) we can rewrite (\ref{gp.18}) as
\begin{equation}
E\(G_{a}G_{b}\prod_{i=1}^{n}G^{2}_{x_{i}}/2\)=\sum_{A\subseteq [1,n]}\mbox{ch}_{A}(a,b)\,\,E\(\prod_{i\notin A}G^{2}_{x_{i}}/2\).\label{gp.20}
\end{equation}
 To see this we use the previous procedure but start with $G_{a}$. Rather than obtain a cycle, since $G_{a}$ appears only once,  eventually we are paired with $G_{b}$. This forms the chain, and the remaining $G^{2}_{x_{i}}$'s lead to cycles as before.
 
 {\it For more details on the material covered in this section, see the beginning of Section 5.1 in \cite{book}. (\ref{gp.12}) is Lemma 5.2.6 in that book, and (\ref{gp.20}) is stated there as (8.93) and proven carefully.}

\section{Markov processes}\label{sec-mp}

Let $S$ be a topological space which is locally compact with countable base (LCCB). Let \[\{p_{t}(x,y), (t,x,y)\in R_{+}^{1}\times S\times S\}\] be a semigroup of sub-probability kernels with respect to some measure $m$ on $S$. That is, 
$p_{t}(x,y)\geq 0$ and satisfies
\begin{equation}
\int p_{t}(x,y)\,dm(y)\leq 1\label{mp.1}
\end{equation}
and
\begin{equation}
\int p_{t}(x,y)p_{s}(y,z)\,dm(y)=p_{t+s}(x,z).\label{mp.1}
\end{equation}
We write $P_{t}$ for the semigroup of operators induced by $p_{t}(x,y)$. 
\begin{equation}
P_{t}f(x)=\int p_{t}(x,y)f(y)\,dm(y),\label{mp.1a}
\end{equation}
and note that  
\begin{equation}
\| P_{t}f \|_{\ff}\leq \| f \|_{\ff}.\label{mp.1b}
\end{equation}

It will be useful to introduce the $\De$ formalism which turns any semigroup of sub-probability kernels $p_{t}(x,y)$ into a semigroup of probability kernels $\wt p_{t}(x,y)$. To do this  we introduce a new point $\De\notin S$, called the cemetery state and extend $m$ to have unit mass at $\De$. Then if we set
 $\wt p_{t}(x,y)= p_{t}(x,y)$ for $y\in S$, 
 $\wt p_{t}(x,\De)=1-\int_{S} p_{t}(x,y)\,dm(y)$,  and $\wt p_{t}(\De,\De)=1$ one can check that the $\wt p_{t}(x,y)$ form a 
semigroup of probability kernels. In the following we will denote by $p_{t}(x,y)$ this extension to a semigroup of probability kernels on $S\cup \De$, and use the convention that for any function $f$ on $S$ we set $f(\De)=0$.

Given such a semigroup of   kernels $p_{t}(x,y)$, we say that  $X=\{X_{t}, t\geq 0\}$ is a  Markov process  with transition densities $p_{t}(x,y)$ if for any 
bounded measurable functions $f_{i}, 1\leq i\leq k$ on $S$, and times $t_{1}<\cdots<t_{k}$
\bea
&&P^{x}\(\prod_{i=1}^{k}f_{i}(X_{t_{i}})\)\label{mp.2}\\
&&=\int p_{t_{1}}(x,y_{1})p_{t_{2}-t_{1}}(y_{1},y_{2})\cdots p_{t_{k}-t_{k-1}}(y_{k-1},y_{k})\prod_{i=1}^{k}f_{i}(y_{i})\,dm(y_{i}).\nn
\eea
Constructing a `nice' Markov process from the kernels $p_{t}(x,y)$ is another story. For now we simply assume that $X$
 has right continuous paths and satisfies the strong Markov property.

For example, for Brownian motion we have $S=R^{1}$, $m$ is Lebesgue measure and  $p_{t}(x,y)=p_{t}(x-y)=e^{-(x-y)^{2}/2t}/\sqrt{2\pi t}$.

We next introduce the $\al$-potential kernels, $\al\geq 0$,
\begin{equation}
u^{\al}(x,y)=\int_{0}^{\ff}e^{-\al t}p_{t}(x, y)\,dt.\label{mp.3}
\end{equation}
We assume that the $u^{\al}(x,y)$ are  continuous for some $\al\geq 0$. 

We note that if $X $ symmetric then   $p_{t}(x, y)$ is positive definite: 
\begin{eqnarray}
\sum_{i,j=1}^{n}a_{i}a_{j}\,\,p_{t}(x_{i}, x_{j})&=& \sum_{i,j=1}^{n}a_{i}a_{j}\,\int \,p_{t/2}(x_{i}, z)\,p_{t/2}(z, x_{j})\,dm(z)
\label{mp.4}\\
&=& \int |\sum_{i=1}^{n}a_{i}\,p_{t/2}(x_{i}, z)|^{2}\,dm(z)\geq 0,  \nonumber
\end{eqnarray}
where the last equality used the  symmetry $p_{t/2}(z, x_{j})=p_{t/2}(x_{j},z)$. This immediately implies that $u^{\al}(x,y)$ is symmetric and positive definite. Hence there exists a  Gaussian process
$G=\{G_{x}, x\in S\}$ with covariance 
\begin{equation}
E\(G_{x} G_{y}\)=u^{\al}(x,y).\label{mp.5}
\end{equation}
Of course, $G$ depends on $\al$.
 When $\al=0$ and $u^{0}$ is finite we refer to $G$
as the Gaussian process associated with $X$.  $G$ is one of the key players in the Isomorphism Theorem. 

We now introduce the other key player, the local time $L=\{L^{y}_{t}, (t,y)\in R_{+}^{1}\times S\}$ defined by
\begin{equation}
L^{y}_{t}=\lim_{\ep\rar 0}\int_{0}^{t}f_{\ep,y}(X_{r})\,dr,\label{mp.6}
\end{equation}
where $f_{\ep,y}$ is an approximate $\de$-function at $y$. That is, $f_{\ep,y}$ is a non-negative function supported in $B(y,\ep)$ with $\int f_{\ep,y}(x)\,dm(x)=1$. If $u^{\al}(x,y)$ is continuous
for some $\al\geq 0$, it can be shown that the limit in (\ref{mp.6}) exists locally uniformly in $t$, $P^{x}$ a.s. It is then easily seen that $L^{y}_{t}$ inherits the following properties from $\int_{0}^{t}f_{\ep,y}(X_{r})\,dr$: $L^{y}_{0}=0$, $L^{y}_{t}$ is continuous and increasing in $t$, and has the additivity property:
\begin{equation}
L^{y}_{t+s}=L^{y}_{t}+L^{y}_{s}\circ\th_{t},\label{mp.7}
\end{equation}
where $\th_{t}\om(r)=\om(r+t)$. 

Thus $L^{y}_{t}$ is continuous in $t$, but what about continuity in $y$? The Isomorphism Theorems allow us to give a complete resolution to this question for symmetric Markov processes.

\subsection{Local time moment formulas}

For  ease of notation we assume that $u^{0}(x,y)$ is continuous, and write it as $u(x,y)$. Our first formula is somewhat similar to  the chain function (\ref{gp.19}) which appears in the Gaussian moment formula (\ref{gp.18}).
\be
 P^{x}\(   \prod_{i=1}^{k}L^{ x_{i}}_{\ff}  \) 
  = \sum_{\pi\in \mathcal{P}_{k}} u(x,x_{\pi(1)})u(x_{\pi(1)},x_{\pi(2)})\cdots   u(x_{\pi(k-1)},x_{\pi(k)}).  \label{mp.12}
 \ee

Proof: It follows from (\ref{mp.2})
\begin{eqnarray}
&&P^{x}\(   \int_{\{0<t_{1}<\cdots<t_{k}<\ff\}}\prod_{i=1}^{k}f_{i}(X_{t_{i}})\,dt_{i}\)
\label{mp.10}\\
&&  =\int u(x,y_{1})u(y_{1},y_{2})\cdots    u(y_{k-1},y_{k}) \prod_{i=1}^{k} f_{i}(y_{i})\,dm (y_{i}). \nonumber
\end{eqnarray}
and consequently, since $R_{+}^{k}=\cup_{\pi\in \mathcal{P}_{k}}\{0<t_{\pi(1)}<t_{\pi(2)}<\cdots<t_{\pi(k)}<\ff\}$ (up to sets of Lebesgue measure $0$),
\begin{eqnarray}
&&P^{x}\(   \prod_{i=1}^{k}\int_{-\ff}^{\ff}f_{i}(X_{t_{i}})\,dt_{i}\)
\label{mp.11}\\
&&  = \sum_{\pi\in \mathcal{P}_{k}}\int u(x,y_{1})u(y_{1},y_{2})\cdots    u(y_{k-1},y_{k}) \prod_{i=1}^{k} f_{\pi(i)}(y_{i})\,dm (y_{i}).\nonumber
\end{eqnarray}
Taking $f_{i}=f_{\ep,x_{i}}$ and then taking the limit as $\ep\to 0$ gives (\ref{mp.12}).

 To prove Dynkin's Isomorphism Theorem we will need a different sort of measure, known as an h-transform of our Markov process $X$. We define a measure $ Q^{x,y}$ by the formula
 \begin{equation}
 Q^{x,y}(F1_{t<\ze} )=P^{x}(F\, u(X_{t}, y)),\hspace{.2 in}  F\in \mathcal{F}_{t}. \label{mp.13}
 \end{equation}
 That is, if we take some functional $F$ which depends only on the path up to time $t$, we first measure $F$ using $P^{x}$, and then, starting at position $X_{t}$, the factor $u(X_{t}, y)$
 measures all possible ways to end up at $y$.
 Here is the moment formula we want:
  \begin{eqnarray}
 &&Q^{x,y}\(   \prod_{i=1}^{k}L^{ x_{i}}_{\ff}  \) 
 \label{mp.14}\\
 &&  = \sum_{\pi\in \mathcal{P}_{k}} u(x,x_{\pi(1)})u(x_{\pi(1)},x_{\pi(2)})\cdots   u(x_{\pi(k-1)},x_{\pi(k)})u(x_{\pi(k)},y).\nonumber
 \end{eqnarray}
 In other words, comparing with (\ref{gp.19}) for the associated process,
   \be
Q^{x,y}\(   \prod_{i=1}^{k}L^{ x_{i}}_{\ff}  \) = \mbox{ch}([1,k]; x,y). \label{mp.15}
 \ee
 
 Proof: If $t_{1}<\cdots<t_{k}$, it follows from the definition (\ref{mp.13}) that
 \bea
&&Q^{x,y}\(\prod_{i=1}^{k}f_{i}(X_{t_{i}})\)=P^{x}(\prod_{i=1}^{k}f_{i}(X_{t_{i}})\, u(X_{t_{k}}, y))\label{mp.2u}\\
&&=\int p_{t_{1}}(x,y_{1})p_{t_{2}-t_{1}}(y_{1},y_{2})\cdots p_{t_{k}-t_{k-1}}(y_{k-1},y_{k})u(y_{k},y)\prod_{i=1}^{k}f_{i}(y_{i})\,dm(y_{i}).\nn
\eea
Hence
\begin{eqnarray}
&&Q^{x,y}\(   \int_{\{0<t_{1}<\cdots<t_{k}<\ff\}}\prod_{i=1}^{k}f_{i}(X_{t_{i}})\,dt_{i}\)
\label{mp.16}\\
&&  =\int u(x,y_{1})u(y_{1},y_{2})\cdots    u(y_{k-1},y_{k})u(y_{k},y)  \prod_{i=1}^{k} f_{i}(y_{i})\,dm (y_{i}). \nonumber
\end{eqnarray}
Arguing as before we then see that
\begin{eqnarray}
&&\hspace{-.3 in}Q^{x,y}\(   \prod_{i=1}^{k}\int_{-\ff}^{\ff}f_{i}(X_{t_{i}})\,dt_{i}\)
\label{mp.17}\\
&&\hspace{-.3 in}  = \sum_{\pi\in \mathcal{P}_{k}}\int u(x,y_{1})u(y_{1},y_{2})\cdots    u(y_{k-1},y_{k})u(y_{k},y) \prod_{i=1}^{k} f_{\pi(i)}(y_{i})\,dm (y_{i}).\nonumber
\end{eqnarray}
Taking $f_{i}=f_{\ep,x_{i}}$ and then taking the limit as $\ep\to 0$ gives (\ref{mp.14}).
 
 {\it For more details about Markov processes and local times, see \cite[Chapter 2]{book}.  (\ref{mp.11}) is Theorem 3.3.2 of that book. The moment formula (\ref{mp.12}) is a special case of Theorem 3.10.1, where we take $T=\ff$, and (\ref{mp.14}) is equivalent to (3.248). }

 \section{The Dynkin Isomorphism Theorem}
 
 The Dynkin Isomorphism Theorem can be expressed as
 \begin{equation}
 E_{G}Q^{x,y}\(F\(   L^{ x_{i}}_{\ff}+{1 \over 2}G^{2}_{ x_{i}} \) \)=  E_{G} \(G_{ x}G_{ y}\, F\(   {1 \over 2}G^{2}_{ x_{i}} \) \).\label{dit.1}
 \end{equation}
 Here, for a bounded measurable function $F$ on $R^{\ff}$ we use the abbreviation
\begin{equation}
F(h_{x_i})=F(h_{x_1},h_{x_2},h_{x_3},\ldots). \label{dit.2}
\end{equation}
$E_{G}$ denotes expectation with respect to  the associated Gaussian process 
$G$.
Note that in (\ref{dit.1}) the associated Gaussian process 
$G$ is independent of the Markov process $X$.   (\ref{dit.1}) is not what is usually referred to as an isomorphism: The right hand side contains only the process $G$, but the left hand side is a mixture of the local time process of $X$ and the independent process $G$. Before giving a proof of (\ref{dit.1}), which will be simple since we have already developed most of tools we need, I would like to give an example to illustrate how to `decouple' $L$ and $G$.

Assume that we know that the associated Gaussian process 
$G$ is a.s. continuous on $S$. We will use the Dynkin Isomorphism Theorem to show that the total local time $L^{z}_{\ff}$ is continuous on $S$, $Q^{x,y}$ a.s. Continuity is a local property, so it is sufficient to show that $L^{z}_{\ff}$ is continuous on any compact subset  $K\subseteq S$, $Q^{x,y}$ a.s. Pick a countable dense subset $D\subseteq K$, and let $F_{D}$ be the indicator function of the event that a function $h$ is uniformly continuous on $D$. Since by assumption $G$ is a.s. continuous on $S$, we have that $F_{D}(G^{2}/2)=1$, a.s. Hence the right hand side of (\ref{dit.1}) is equal to $ E_{G} \(  G_{ x}G_{ y}\)=u(x,y)$, which is precisely the total mass of the measure $ E_{G}Q^{x,y}$. Therefore $F_{D}(L_{\ff}+G^{2}/2)=1$, a.s. That is, $L^{z}_{\ff}+G_{z}^{2}/2$ is almost surely uniformly continuous on $D$, and since we know this is true of $G_{z}^{2}/2$ we have established that $L^{z}_{\ff}$ is almost surely uniformly continuous on $D$. This is basically what we wanted to show. Standard techniques allow us to extend $L^{z}_{\ff}$ by continuity to $K$, and verify that this extension is indeed the total local time $L^{z}_{\ff}, z\in K$.

By the way, this result is not purely academic. Necessary and sufficient conditions for the a.s. continuity of a Gaussian process in terms of its covariance are known. We describe this in the next section.

Proof of the Dynkin Isomorphism Theorem: We first take $F$ to be a product, and show that
 \begin{equation}
 E_{G}Q^{x,y}\(\prod_{i=1}^{k}\(   L^{ x_{i}}_{\ff}+{1 \over 2}G^{2}_{ x_{i}} \) \)=  E_{G} \(G_{ x}G_{ y} \prod_{i=1}^{k}   {1 \over 2}G^{2}_{ x_{i}} \).\label{dit.3}
 \end{equation}
 Expanding the product on the left hand side, (\ref{dit.3}) is 
  \begin{equation}
\sum_{A\subseteq [1,k]} Q^{x,y}\(\prod_{i\in A}   L^{ x_{i}}_{\ff}   \)E_{G}\(\prod_{i\notin A}{1 \over 2}G^{2}_{ x_{i}}\)=  E_{G} \(G_{ x}G_{ y}\, \prod_{i=1}^{k}   {1 \over 2}G^{2}_{ x_{i}}  \).\label{dit.3}
 \end{equation}
In view of (\ref{mp.15}) this is just (\ref{gp.20}).

To extend this to general bounded measurable $F$, we need only show that the two sides of (\ref{dit.1}) are determined by their moments and this will follow once we show that both $L^{z}_{\ff}$ and $G_{z}^{2}$ are exponentially integrable. But it follows from (\ref{mp.14}) that
\begin{equation}
Q^{x,y}\(\(L^{z}_{\ff}\)^{n}\)=n! u(x,z)(u(z,z))^{n-1}u(z,y)\label{dit.4}
\end{equation}
and from (\ref{gp.12}) that
\begin{equation}
\label{dit.5}
E\(G^{2n}_{z}\)=|\mathcal{R}_{2n}|u^{n}(z,z),
\end{equation}
and $|\mathcal{R}_{2n}|$, the number of pairings of $2n$ objects, is bounded by $n!c^{n}$. (In fact, the exponential integrability of the square of a normal random variable is well know and easy to compute explicitly.)
\qed

{\it     The proof of  the Dynkin Isomorphism Theorem given here is found in \cite[Section 8.3.1]{book}.    }

 \section{The Eisenbaum Isomorphism Theorem}

One problem with  the Dynkin Isomorphism Theorem is the appearance of the measure $Q^{x,y}$.
The following  Isomorphism Theorem of Eisenbaum deals with the natural measure $P^{x}$, but at some cost. It says that for any $s>0$
 \begin{equation}
 E_{G}P^{x}\(F\(   L^{ x_{i}}_{\ff}+{1 \over 2}(G_{ x_{i}}+s)^{2} \) \)=  E_{G} \(\(1+{G_{x} \over s}\) F\(   {1 \over 2}(G_{ x_{i}}+s)^{2} \) \).\label{dit.6}
 \end{equation}

Proof: Once again it suffices to prove this when $F$ is a product, in which case, after expanding the first factor on the right hand side of (\ref{dit.6}), it takes the form
  \bea
&&\sum_{A\cup B =[1,k]} P^{x}\(\prod_{i\in A}   L^{ x_{i}}_{\ff}   \)E_{G}\(\prod_{i\in B}{1 \over 2}(G_{ x_{i}}+s)^{2}\)\label{dit.7}\\
&&=  E_{G} \( \prod_{i=1}^{k}\(   {1 \over 2}(G_{ x_{i}}+s)^{2} \) \)+E_{G} \( {G_{x} \over s} \prod_{i=1}^{k}\(   {1 \over 2}(G_{ x_{i}}+s)^{2} \) \).\nn
 \eea
The first term on the right hand side corresponds to the term on the on the left hand side
with $A=\emptyset$. If $A\not=\emptyset$, recall that by (\ref{mp.12}), if $A=\{a_{1}, a_{2},\ldots,    a_{|A|}\}$
\be
 P^{x}\(   \prod_{i\in A}L^{ x_{i}}_{\ff}  \) 
  = \sum_{\pi\in \mathcal{P}_{|A|}} u(x,x_{a_{\pi(1)}})u(x_{a_{\pi(1)}},x_{a_{\pi(2)}})\cdots   u(x_{a_{\pi(|A|-1)}},x_{a_{\pi(|A|)}}).  \label{dit.8}
 \ee
 For the expectation on the right of (\ref{dit.7}), start with ${G_{x} \over s}$ and apply the Gaussian moment formula (\ref{gp.12}). $G_{x}$ must be paired with something. It can be paired with one of the two factors of $G_{x_{a_{\pi(1)}}}$, canceling the $1/2$ and giving rise to the factor $u(x,x_{a_{\pi(1)}})$. The other factor $G_{x_{a_{\pi(1)}}}$ might be paired with one of the two factors of $G_{x_{a_{\pi(2)}}}$, canceling the $1/2$ and giving rise to the factor $u(x_{a_{\pi(1)}},x_{a_{\pi(2)}})$.
We proceed in the way until we pair $G_{x_{a_{\pi(|A|-1)}}}$ with $G_{x_{a_{\pi(|A|)}}}$ from  one of the two factors of 
$(G_{ x_{a_{\pi(|A|)}}}+s)^{2}$. From the other factor we take $s$, canceling the $1/s$ from 
${G_{x} \over s}$. Thus we have obtained (\ref{dit.8}) and what remains from this expectation on the right of (\ref{dit.7}) is precisely $E_{G}\(\prod_{i\in B}{1 \over 2}(G_{ x_{i}}+s)^{2}\)$.\qed

This completes the proof of the Eisenbaum Isomorphism Theorem, but it is of interest, and will be useful later on, to figure out explicitly the other terms. We show that
\begin{equation}
E_{G} \( \prod_{i=1}^{k}\(   {1 \over 2}(G_{ x_{i}}+s)^{2} \) \)=\sum_{\stackrel{A_{1}\cup\cdots\cup A_{l}} {\cup B_{1}\cup\cdots\cup B_{m}=[1,k]}}
\prod_{i=1}^{l}{1 \over 2} \mbox{cy} (A_{i})    \prod_{j=1}^{m}{s^{2} \over 2}\mbox{ch} (B_{j})
\label{dit.9}
\end{equation}
where $\mbox{cy} (A)$ is defined in (\ref{gp.16}), $\mbox{ch}(B)=1$ if $|B|=1$ and, if  $|B|>1$ with $B=\{b_{1},b_{2},\cdots, b_{|B|}\}$  then the chain function $\mbox{ch}(B)$ is defined as 
\begin{equation}
\mbox{ch}(B)=\sum_{\pi\in \mathcal{P}_{| B|}}u(x_{b_{\pi(1)}},  x_{b_{\pi(2)}})\cdots u(x_{b_{\pi (| B|-1)}},  x_{b_{\pi (| B|)}}).\label{dit.10}
\end{equation} 
Note that the `chains' in $\mbox{ch}(B)$ are oriented. For example if $B=\{1,2\}$ then 
$\mbox{ch}(B)=u(x_{1}, x_{2})+u(x_{2}, x_{1})$. For the symmetric case we are dealing with this is
$2u(x_{1}, x_{2})$.

Proof of (\ref{dit.9}): It will be convenient to rewrite this as
\begin{eqnarray}
&&E_{G} \( \prod_{i=1}^{k}\(   {1 \over 2}(G_{ x_{i}}+s)^{2} \) \)
\label{dit.11}\\
&&=\sum_{A\cup B =[1,k]}\( \sum_{ A_{1}\cup\cdots\cup A_{l}=A}\,
\prod_{i=1}^{l}{1 \over 2} \mbox{cy} (A_{i})  \sum_{B_{1}\cup\cdots\cup B_{m}=B} \, \prod_{j=1}^{m}{s^{2} \over 2}\mbox{ch} (B_{j})\).   \nonumber
\end{eqnarray}
There are many terms in the expansion of $\prod_{i=1}^{k}\(   {1 \over 2}(G_{ x_{i}}+s)^{2} \)$. If we look at $\prod_{i\in A}   {1 \over 2}G_{ x_{i}}^{2} $ and pair together  all factors in this product, then using (\ref{gp.17}) we obtain the term on the right hand side of (\ref{dit.11}) containing cycles. To obtain the term involving chains, if for example $B_{j}=\{b_{1},b_{2},\cdots, b_{|B_{j}|}\}$ we can obtain $u(x_{b_{\pi(1)}},  x_{b_{\pi(2)}})\cdots u(x_{b_{\pi (| B|-1)}},  x_{b_{\pi (| B|)}})$ by looking at a specific pairing of $H=:sG_{x_{b_{\pi(1)}}}  \(\prod_{i=2}^{|B_{j}|-1}  {1 \over 2}   G^{2}_{x_{b_{\pi(i)}}}\)       sG_{x_{b_{\pi(|B_{j}|)}}} $. That is, we pair $G_{x_{b_{\pi(1)}}}$
with one of the two factors $G_{x_{b_{\pi(2)}}}$, pair the other factor $G_{x_{b_{\pi(2)}}}$
with one of the two factors $G_{x_{b_{\pi(3)}}}$, until finally we pair the remaining factor $G_{x_{b_{\pi(|B_{j}|-1)}}} $ with $G_{x_{b_{\pi(|B_{j}|)}}} $. In this way we have cancelled all the factors of 
$1/2$ in $H$ and obtained a factor of $s^{2}$. But note that this pairing is unoriented, while as we mentioned the `chains' in $\mbox{ch}(B_{j})$ are oriented. This accounts for the factor $1/2$
multiplying $\mbox{ch}(B_{j})$ in (\ref{dit.11}).\qed

\subsection{Bounded discontinuities}

We now present an application of the Eisenbaum Isomorphism Theorem. We first recall the fundamental result of Talagrand that a Gaussian process
$G= \{G_{x}, x\in S\}  $
is continuous a.s. if and only if there exists a probability measure $\nu$ on $S$ such that
\begin{equation}
\lim_{\de\rar 0}\sup_{s\in S}\int_{0}^{\de}\(\log {1 \over \nu (B_{d}(s,u))}\)^{1/2}du=0.\label{dit.12}
\end{equation} 
Here, continuity is with respect to the metric 
$d(x,y)=\(E\(\lc G_{x}- G_{y}\rc^{2}\)\)^{1/2}$ which can be expressed in terms of the covariance 
$u(x,y)$ of $G$.
 
  Marcus and I  used  Isomorphism Theorems to show  that for symmetric Markov processes with continuous potential densities, the total local time $L=\{L^{z}_{\ff}, z\in S\}$ will be $P^{x}$ almost surely continuous for each $x\in S$ if and only if the associated Gaussian process $G$
is almost surely continuous. By the result of Talagrand  we have an explicit condition in terms of 
 the potential densities $u(x,y)$.

We have already indicated how to use  Isomorphism Theorems to show that if the associated Gaussian process $G$
is almost surely continuous, then the total local time $L=\{L^{z}_{\ff}, z\in S\}$ will be almost surely continuous. We now show how to use the Eisenbaum Isomorphism Theorem to go in the other direction, that is, to show that  if the associated Gaussian process $G$
is not almost surely continuous, then the total local time $L=\{L^{z}_{\ff}, z\in S\}$ cannot be  be $P^{x}$ almost surely continuous for each $x\in S$. The key to this result is the fact that a 
Gaussian process can only be discontinuous in very special ways, which we now recall.

Set
\begin{equation}
M_{f}(x_{0})=\lim_{\ep\to 0}\sup_{x\in B_{d}(x_{0},\ep)}f(x),\hspace{.2 in}m_{f}(x_{0})=\lim_{\ep\to 0}\inf_{x\in B_{d}(x_{0},\ep)}f(x).\label{dit.13}
\end{equation}
Let $G= \{G_{x}, x\in S\}  $ be a Gaussian process with continuous covariance. If $G$ is not almost surely continuous then there exists $x_{0}\in S$, a $\bb(x_{0})>0$ and a countable dense subset $C\subseteq S$ such that
\begin{equation}
M_{G_{|C}}(x_{0})=G_{x_{0}}+\bb(x_{0})\hspace{.2 in}\mbox{and}\hspace{.2 in}m_{G_{|C}}(x_{0})=G_{x_{0}}-\bb(x_{0})\hspace{.2 in}a.s.\label{dit.14}
\end{equation}

When $0<\bb(x_{0})<\ff$ we say that $G$ has a bounded discontinuity at $x_{0}$. If 
$\bb(x_{0})=\ff$ we say that $G$ has an unbounded discontinuity at $x_{0}$. We will now use the Eisenbaum Isomorphism Theorem to show that if $G$ has a bounded discontinuity at $x_{0}$
then $L=\{L^{z}_{\ff}, z\in S\}$ will be discontinuous at $x_{0}$, $P^{x_{0}}$ almost surely. The case of an unbounded discontinuity is somewhat more complicated and we refer the interested reader to \cite[Chapter 9.2]{book}.

Proof: Simple algebra shows that
\begin{equation}
(G_{ x}+s)^{2} -(G_{ x_{0}}+s)^{2} =(G_{ x}-G_{ x_{0}})^{2}+2(G_{ x_{0}}+s)(G_{ x}-G_{ x_{0}}). \label{dit.15}
\end{equation}
Using this we claim that almost surely
\begin{equation}
\lim_{\ep\to 0}\sup_{x\in B_{d}(x_{0},\ep)\cap C}(G_{ x}+s)^{2} -(G_{ x_{0}}+s)^{2} =\bb^{2}(x_{0})+2\bb(x_{0})|G_{ x_{0}}+s|.\label{dit.16a}
\end{equation}
To see this, look at the right hand side of (\ref{dit.15}). If $G_{ x_{0}}+s>0$ we obtain (\ref{dit.16a}) by taking a sequence of points $x_{n}\rar x_{0}$ such that, by (\ref{dit.14}), $G_{ x_{n}}-G_{ x_{0}}\rar \bb(x_{0})$, while if $G_{ x_{0}}+s<0$ we obtain (\ref{dit.16a}) by taking a sequence of points $x_{n}\rar x_{0}$ such that, by (\ref{dit.14}), $G_{ x_{n}}-G_{ x_{0}}\rar -\bb(x_{0})$.

We can rewrite (\ref{dit.16a}) in form more appropriate to the Eisenbaum Isomorphism Theorem: 
\begin{equation}
\lim_{\ep\to 0}\sup_{x\in B_{d}(x_{0},\ep)\cap C}{1 \over 2}(G_{ x}+s)^{2} -{1 \over 2}(G_{ x_{0}}+s)^{2} ={\bb^{2}(x_{0})\over 2}+2^{1/2}\bb(x_{0})\sqrt{{1 \over 2}(G_{ x_{0}}+s)^{2}},\label{dit.16}
\end{equation}
almost surely. Let $x_{i}$ be an enumeration of the points in $C$. We now apply the Eisenbaum Isomorphism Theorem with $F(f_{x_{i}})$ the indicator function of the event  
\begin{equation}
\lim_{\ep\to 0}\sup_{x_{i}\in B_{d}(x_{0},\ep)}f_{x_{i}} -f_{x_{0}} ={\bb^{2}(x_{0})\over 2}+2^{1/2}\bb(x_{0})\sqrt{f_{x_{0}}}.\label{dit.17}
\end{equation}
By (\ref{dit.16}), $F({1 \over 2}(G_{ x_{i}}+s)^{2})=1$ a.s. The Eisenbaum Isomorphism Theorem
then implies that for any $s>0$, $F(L^{x_{i}}_{\ff}+{1 \over 2}(G_{ x_{i}}+s)^{2})=1$ a.s. That is,
\begin{eqnarray}
&&\lim_{\ep\to 0}\sup_{x_{i}\in B_{d}(x_{0},\ep)} \(L^{x_{i}}_{\ff}-L^{x_{0}}_{\ff}+{1 \over 2}(G_{ x_{i}}+s)^{2} -{1 \over 2}(G_{ x_{0}}+s)^{2}\)
\label{dit.18}\\
&&={\bb^{2}(x_{0})\over 2}+2^{1/2}\bb(x_{0})\sqrt{L^{x_{0}}_{\ff}+{1 \over 2}(G_{ x_{0}}+s)^{2}}.   \nonumber\\
&&\geq {\bb^{2}(x_{0})\over 2}+2^{1/2}\bb(x_{0})\sqrt{L^{x_{0}}_{\ff}}.   \nonumber
\end{eqnarray}
Using the fact that $\lim\sup_{i}A_{i}+\lim\sup_{i}B_{i}\geq \lim\sup_{i}(A_{i}+B_{i})$ and then (\ref{dit.16}) we see that almost surely
\begin{equation}
\lim_{\ep\to 0}\sup_{x_{i}\in B_{d}(x_{0},\ep)} L^{x_{i}}_{\ff}-L^{x_{0}}_{\ff}\geq 2^{1/2}\bb(x_{0})\sqrt{L^{x_{0}}_{\ff}}-2^{1/2}\bb(x_{0})|G_{ x_{0}}+s|.\label{dit.19}
\end{equation}

At first glance this doesn't seem very useful. We want to show that $L$ has a discontinuity at $x_{0}$, that is, that the left hand side is strictly positive, but because we are subtracting $2^{1/2}\bb(x_{0})|G_{ x_{0}}+s|$ on the right hand side, the right hand side might be negative! 

I will now perform a magic trick. I will make the $\bb(x_{0})|G_{ x_{0}}+s|$ disappear before your very eyes! For this purpose recall that $L$ and $G$ are independent and in fact live on different spaces. To emphasize this we write (\ref{dit.19}) as the statement that
\begin{equation}
\lim_{\ep\to 0}\sup_{x_{i}\in B_{d}(x_{0},\ep)} L^{x_{i}}_{\ff}(\om)-L^{x_{0}}_{\ff}(\om)\geq 2^{1/2}\bb(x_{0})\sqrt{L^{x_{0}}_{\ff}(\om)}-2^{1/2}\bb(x_{0})|G_{ x_{0}}(\om')+s|\label{dit.19}
\end{equation}
almost surely with respect to $P^{x}\times E_{G}$. By Fubini's theorem then, this holds $P^{x}$
almost surely for $P_{G}$ almost every $\om'$. That is, (\ref{dit.19}) holds $P^{x}$
almost surely for all $\om'\in \Om'$ where $P_{G}( \Om')=1$. If we could only find an $\om'\in \Om'$ with $|G_{ x_{0}}(\om')+s|=0$ we would be done. We do something similar. Fix $\de>0$ and set $s=\de$. Since $G_{ x_{0}}$ is a normal random variable, we have
\begin{equation}
P_{G}( |G_{ x_{0}}|\leq \de)>0.\label{dit.20}
\end{equation}
Since $P_{G}( \Om')=1$ we can find $\om'\in \Om'$ with $ |G_{ x_{0}}(\om')|\leq \de$. By the above we now see that  $P^{x}$
almost surely
\begin{equation}
\lim_{\ep\to 0}\sup_{x_{i}\in B_{d}(x_{0},\ep)} L^{x_{i}}_{\ff}(\om)-L^{x_{0}}_{\ff}(\om)\geq 2^{1/2}\bb(x_{0})\sqrt{L^{x_{0}}_{\ff}(\om)}-2^{1/2}\bb(x_{0})2\de. \label{dit.21}
\end{equation}
Since this is true for any $\de>0$ we have in fact shown that $P^{x}$
almost surely
\begin{equation}
\lim_{\ep\to 0}\sup_{x_{i}\in B_{d}(x_{0},\ep)} L^{x_{i}}_{\ff}(\om)-L^{x_{0}}_{\ff}(\om)\geq 2^{1/2}\bb(x_{0})\sqrt{L^{x_{0}}_{\ff}(\om)}. \label{dit.21}
\end{equation}
Is $L^{x_{0}}_{\ff}(\om)>0$ almost surely? This depends on whether or not the path has visited $x_{0}$. But we can give a simple proof that $L^{x_{0}}_{\ff}>0$, $P^{x_{0}}$ almost surely, so by the above we can conclude that the local time $L$ is discontinuous $P^{x_{0}}$ almost surely. This is all that we wanted to establish.

The fact   that $L^{x_{0}}_{\ff}>0$, $P^{x_{0}}$ almost surely follows from (\ref{mp.12}) which implies that for all $n$
\begin{equation}
P^{x_{0}}\(\(L^{x_{0}}_{\ff}\)^{n}\)=n!(u(x_{0},x_{0}))^{n},\label{dit.22}
\end{equation}
which implies that  $L^{x_{0}}_{\ff}$ is distributed  under $P^{x_{0}}$ as an exponential random variable (with mean $u(x_{0},x_{0})$). Since exponential random variables are strictly positive almost surely, we are done. \qed

{\it     The proof of  the Eisenbaum Isomorphism Theorem given here is similar to that   in \cite[Section 8.3.2]{book}.  For  Talagrand's theorem, see Chapter 6 of that book. The property (\ref{dit.14}) concerning discontinuities of Gaussian processes is Theorem 5.3.7 and our proof that if the associated Gaussian process has a bounded discontinuity, the local time will be discontinuous is given in detail in Chapter 9.1 of the book. }

\section{The generalized second Ray-Knight theorem}

We fix some point in $S$ which we denote by $0$. Set
\begin{equation}
\tau(t)=\inf \{s\,|\,L^{ 0}_{s}>t\},\label{rk.1}
\end{equation}
the right continuous inverse local time at $0$, and
\begin{equation}
T_{y}=\inf \{s\,|\,X_{s}=y\},\label{rk.0}
\end{equation}
the first hitting time of $y$.
 For this section we assume that $u^{\al}(x,y)$ is continuous for any $\al>0$ and  $u(0,0)=\ff$. In addition, we assume that  $P^{0}\(T_{x}<\ff\)=P^{x}\(T_{0}<\ff\)=1$ for all $x\in S$.
The generalized second  Ray-Knight theorem states that for all $t>0$
\begin{equation}
E_{G}P^{0}\(F\(   L^{ x_{i}}_{\tau(t)}+{1 \over 2}\eta^{2}_{ x_{i}} \) \)= E_{G}\(F\( {1 \over 2}(\eta_{ x_{i}}+\sqrt{2t})^{2} \) \),\label{rk.2}
\end{equation}
where $\{\eta_{x}, x\in S\}$ is the Gaussian process with covariance
\begin{equation}
u_{T_{0}}(x,y)=E^{x}\(L^{ y}_{T_{0}}\).\label{rk.3}
\end{equation}
 We will see below that indeed $u_{T_{0}}(x,y)$ is symmetric and positive definite. The notation 
 $u_{T_{0}}(x,y)$ is meant to suggest that $u_{T_{0}}(x,y)$
  is the potential density of the Markov process obtained by killing $X$ at $T_{0}$. This is true, \cite[Chapter 4.5]{book}, but we will not use this fact.
 
 Recall that by  (\ref{dit.9}), after replacing $s$ by $\sqrt{2t}$.
 \begin{equation}
E_{G} \( \prod_{i=1}^{k}\(   {1 \over 2}(\eta_{ x_{i}}+\sqrt{2t})^{2}\) \)=\sum_{\stackrel{A_{1}\cup\cdots\cup A_{l}} {\cup B_{1}\cup\cdots\cup B_{m}=[1,k]}}
\prod_{i=1}^{l}{1 \over 2} \mbox{cy}_{0} (A_{i})    \prod_{j=1}^{m}t \mbox{ ch}_{0} (B_{j}).
\label{rk.4}
\end{equation}
Here $\mbox{cy}_{0}$ and $\mbox{ ch}_{0}$ refer to the covariance $u_{T_{0}}(x,y)$ of $\eta$. We emphasize that the sum is over unordered partitions of $[1,k]$. That is, $B_{1}=\{1,2,3\},B_{2}=\{4,5\}$ and $B_{2}=\{1,2,3\},B_{1}=\{4,5\}$  are not counted separately in the sum.
Using (\ref{gp.17}) as before, the generalized second  Ray-Knight theorem will be proven once we show that for all $t>0$
\be
 P^{0}\(   \prod_{i=1}^{k}L^{ x_{i}}_{\tau(t)}  \) 
  = \sum_{m=1}^{k}\sum_{\stackrel{\mbox{\scriptsize unordered}}{B_{1}\cup\cdots\cup B_{m}=[1,k]}} \,\,t^{m}\prod_{j=1}^{m} \mbox{ ch}_{0} (B_{j}).  \label{rk.5}
 \ee 
This, however, is not so simple.   Our local time moment formulas are for the total local time of a Markov process, but $X$ killed at $\tau(t)$ is not a Markov process. To prove (\ref{rk.5}) we let $\la$
be an independent exponential random variable with mean $\al$ and show that
\be
 P^{x}_{\la}\(   \prod_{i=1}^{k}L^{ x_{i}}_{\tau(\la)}  \) 
  = \sum_{\pi\in \mathcal{P}_{k}} u_{\tau(\la)}(x,x_{\pi(1)})u_{\tau(\la)}(x_{\pi(1)},x_{\pi(2)})\cdots   u_{\tau(\la)}(x_{\pi(k-1)},x_{\pi(k)})  \label{rk.6}
 \ee
 where $ P^{x}_{\la}= P^{x}\times  P_{\la}$ and 
 \begin{equation}
u_{\tau(\la)}(x,y):=u_{T_{0}}(x,y)+\al.\label{rk.9}
\end{equation}
Once again, the notation 
 $u_{\tau(\la)}(x,y)$ is meant to suggest that $u_{\tau(\la)}(x,y)$
  is the potential density of a symmetric Markov process with probabilities $ P^{x}_{\la}$ obtained by killing $X$ at $\tau(\la)$. And once again this   is true, \cite{five},  but we will give a proof of  (\ref{rk.6})-(\ref{rk.9})   which does not use this fact.
  
Combining (\ref{rk.6}) and (\ref{rk.9}) and expanding the product we see that
\begin{eqnarray}
&& P^{0}_{\la}\(   \prod_{i=1}^{k}L^{ x_{i}}_{\tau(\la)}  \) 
\label{rk.10}\\
&&  = \sum_{\pi\in \mathcal{P}_{k}}\al \(u_{T_{0}}(x_{\pi(1)},x_{\pi(2)})+\al\)\cdots \(u_{T_{0}}(x_{\pi(k-1)},x_{\pi(k)})+\al\)  \nonumber\\
&&  = \sum_{m=1}^{k} \sum_{\stackrel{\mbox{\scriptsize ordered}}{B_{1}\cup\cdots\cup B_{m}=[1,k]}}\,\,\prod_{j=1}^{m}\al \mbox{ ch}_{0} (B_{j}),  \nonumber
\end{eqnarray}
where now the sum is over ordered partitions, since it comes from a sum over permutations, where order counts. Thus we can write
\begin{equation}
 P^{0}_{\la}\(   \prod_{i=1}^{k}L^{ x_{i}}_{\tau(\la)}  \) =\sum_{m=1}^{k} \sum_{\stackrel{\mbox{\scriptsize unordered}}{B_{1}\cup\cdots\cup B_{m}=[1,k]}}\,\,\al ^{m}m!\prod_{j=1}^{m}\mbox{ ch}_{0} (B_{j}).\label{rk.11}
\end{equation}
Since $\int_{0}^{ \ff}e^{ -t/\al}t^{ m}\,dt/\al=\al^{ m}m!$, we have shown that for any $\al>0$
\bea
&&
\int_{0}^{\ff} e^{-t/\al} P^{0}\(   \prod_{i=1}^{k}L^{ x_{i}}_{\tau(t)}  \) \,dt/\al  \label{rk.12}\\
&&
  =\int_{0}^{\ff} e^{-t/\al}  \sum_{m=1}^{k}\sum_{\stackrel{\mbox{\scriptsize unordered}}{B_{1}\cup\cdots\cup B_{m}=[1,k]}} \,\,t^{m}\prod_{j=1}^{m} \mbox{ ch}_{0} (B_{j})\,dt/\al. \nn
 \eea 
Since $\tau(t)$ is right continuous, we have established (\ref{rk.5}) and hence the generalized second  Ray-Knight theorem, (\ref{rk.2}).

We still have to fill in some missing pieces. We first show that $u_{T_{0}}(x,y)$ is symmetric and positive definite. In fact, we show that
\begin{equation}
u_{T_{0}}(x,y)=\lim_{\al\to 0}\(u^{\al}(x,y)-{u^{\al}(x,0)u^{\al}(0,y) \over u^{\al}(0,0)}\).\label{rk.13}
\end{equation}
This will show that $u_{T_{0}}(x,y)$ is symmetric, and since if $G_{x}$ is the Gaussian process with covariance $u^{\al}(x,y)$ then
\begin{equation}
E\(\(G_{x}-{u^{\al}(x,0) \over u^{\al}(0,0)}G_{0}\)\(G_{y}-{u^{\al}(y,0) \over u^{\al}(0,0)}G_{0}\)\)=u^{\al}(x,y)-{u^{\al}(x,0)u^{\al}(0,y) \over u^{\al}(0,0)},\label{rk.13a}
\end{equation}
(\ref{rk.13}) will also show that $u_{T_{0}}(x,y)$ is positive definite.

Before showing (\ref{rk.13}), let us illustrate it for Brownian motion. 
We first show that
\begin{equation}
u^{\al}(x,y)={e^{-\sqrt{2\al}\,|x-y|} \over \sqrt{2\al}}.\label{rk.14}
\end{equation}
To see this we use the Fourier representation
\bea
&&
u^{\al}(x,y)=\int_{0}^{\ff}e^{-\al t}p_{t}(x,y)\,dt\label{rk.15}\\
&&={1 \over 2\pi}\int_{0}^{\ff}e^{-\al t}\(\int_{-\ff}^{\ff}e^{iz  (x-y)-tz^{2}/2} \,dz\)\,dt={1 \over 2\pi}\int_{-\ff}^{\ff} {e^{iz  (x-y)} \over \al+z^{2}/2}\,dz.
\nn
\eea
$ \al+z^{2}/2$ has two roots in the complex plane $z_{\pm}=\pm i\sqrt{2\al}$. When $(x-y)>0$ we can  evaluate the right hand term in (\ref{rk.15}) by using the residue at $z_{+}$, while if $(x-y)<0$
we use the residue at $z_{-}$. This proves (\ref{rk.14}). Applying this to (\ref{rk.13}) we see that
\begin{eqnarray}
u_{T_{0}}(x,y)&=&\lim_{\al\to 0}\({e^{-\sqrt{2\al}\,|x-y|} \over \sqrt{2\al}}-{e^{-\sqrt{2\al}\,(|x|+|y|)} \over \sqrt{2\al}}\)
\label{rk.16}\\
&=&  (|x|+|y|)-|x-y|=2\,\,( |x|\wedge |y|).\nonumber
\end{eqnarray}
Thus the process $\{\eta_{x}, x\in R^{1}\}$ corresponding to Brownian motion is just $\sqrt{2}$ times two-sided Brownian motion.

We now return to the proof of (\ref{rk.13}). Let $Z$ be the process obtained by killing $X$
at an independent exponential time $\rho$ of mean $1/\bb$. That is, $Z_{t}(\om)=X_{t}(\om)$ if $t<\la$ and $Z_{t}(\om)=\De$ if $t\geq \rho$. Then, recalling our convention that for a function on $S$ we take $f(\De)=0$,
\begin{equation}
E^{x}_{\rho}\(f(Z_{t})\)=E^{x}_{\rho}\(1_{\{\rho>t\}}f(X_{t})\)=e^{-\bb t}E^{x}\(f(X_{t})\)=e^{-\bb t}\int p_{t}(x,y)f( y)\,dm( y).\label{rk.17}
\end{equation}
It follows that $Z$ is a symmetric Markov process whose $0$-potential density is $u^{\bb }(x,y)$. And since the total local time of $Z$ at $y$ is $L^{y}_{\rho}$, see (\ref{mp.6}), we have
\be
u^{\bb }(x,y)=E^{x}_{\rho}\(L^{y}_{\rho}\)=\int_{0}^{\ff}e^{-\bb t}E^{x}\(L^{y}_{t}\)\,\bb \,dt
=E^{x}\( \int_{0}^{\ff}e^{-\bb t}L^{y}_{t}\,\bb \,dt\).\label{rk.18}
\ee
Since $L^{y}_{t}$ is continuous and increasing in $t$, integration by parts then shows that
\be
u^{\bb}(x,y)
=E^{x}\( \int_{0}^{\ff}e^{-\bb t}\,dL^{y}_{t}\).\label{rk.19}
\ee

Using again the additivity of local time and then the strong Markov property we see that for any $\al>0$
\begin{eqnarray}
E^{x}\( \int_{T_{0}}^{\ff}e^{-\al t}\,dL^{y}_{t}\)&=&E^{x}\(e^{-\al T_{0}} \(\int_{0}^{\ff}e^{-\al t}\,dL^{y}_{t}\)\circ\th_{T_{0}}\)
\label{rk.20}\\
&=&  E^{x}\(e^{-\al T_{0}} \)u^{\al}(0,y), \nonumber
\end{eqnarray}
where we use the fact that  $P^{x}(T_0<\ff)=1 $ and that  $X_{ T_0}=0$.
In particular, for $y=0$, recalling that  $L^{0}_{t}$ cannot grow until time $T_{0}$,  this gives  
\begin{equation}
u^{\al}(x,0)=E^{x}\( \int_{T_{0}}^{\ff}e^{-\al t}\,dL^{0}_{t}\)=E^{x}\(e^{-\al T_{0}} \)u^{\al}(0,0),\label{rk.21}
\end{equation}
showing that $E^{x}\(e^{-\al T_{0}} \)=u^{\al}(x,0)/u^{\al}(0,0)$. Putting this back into (\ref{rk.20})
we obtain
\begin{equation}
E^{x}\( \int_{T_{0}}^{\ff}e^{-\al t}\,dL^{y}_{t}\)={u^{\al}(x,0)u^{\al}(0,y) \over u^{\al}(0,0)}.\label{rk.22}
\end{equation}
Together with (\ref{rk.19}) this shows that
\begin{equation}
E^{x}\( \int_{0}^{T_{0}}e^{-\al t}\,dL^{y}_{t}\)=u^{\al}(x,y)-{u^{\al}(x,0)u^{\al}(0,y) \over u^{\al}(0,0)},\label{rk.23}
\end{equation}
and letting $\al\to 0$ completes the proof of (\ref{rk.13}).

The proof of (\ref{rk.6})-(\ref{rk.9}) is more complicated and we defer the proof until after we present an application of the generalized second Ray-Knight theorem. However, we point out that 
if we knew that the process obtained by killing $X$ at $\tau(\la)$ is a symmetric Markov process with continuous potential densities $u_{\tau(\la)}(x,y)$, then (\ref{rk.6}) would simply be our moment formula (\ref{mp.12}), and in particular we would have
\begin{equation}
E^{x}_{\la}\(L^{y}_{\tau(\la)}\)=u_{\tau(\la)}(x,y).\label{rk.6n}
\end{equation}
Since $\tau(t)$ cannot grow until the process first reaches $0$ we have $\tau(\la)=T_{0}+\tau(\la)\circ \th_{T_{0}}$. Hence, using (\ref{mp.7}), the additivity of local times,
\begin{eqnarray}
 E^{x}_{\la}\(L^{y}_{\tau(\la)}\)
\label{rk.7} 
&=&E^{x}_{\la}\(L^{y}_{T_{0}+\tau(\la)\circ \th_{T_{0}}}\)  \nonumber\\
&=&  E^{x}\(L^{y}_{T_{0}}\)+ E^{x}_{\la}\(L^{y}_{\tau(\la)\circ \th_{T_{0}}}\circ \th_{T_{0}}\)\nn\\
&=&  E^{x}\(L^{y}_{T_{0}}\)+ E^{0}_{\la}\(L^{y}_{\tau(\la)}\),\nn
\end{eqnarray}
where the last step used the strong Markov property at the stopping time $T_{0}$. Thus 
\[u_{\tau(\la)}(x,y)=u_{T_{0}}(x,y)+u_{\tau(\la)}(0,y).\] By symmetry we   see that 
\begin{equation}
u_{\tau(\la)}(0,y)=u_{\tau(\la)}(y,0)=E^{y}_{\la}\(L^{0}_{\tau(\la)}\)=E_{\la}\(\la\)=\al.\label{rk.8}
\end{equation}
which combined with the previous display  gives (\ref{rk.9}).

\subsection{Favorite points}

We now illustrate how one can apply  the generalized second Ray-Knight theorem. Let $X$
be a Markov process in $R^{1}$ with continuous potential densities $u^{\al}(x,y)$ and jointly continuous local times $L^{ x}_{t}$. Let
\begin{equation}
\mathcal{V}_{t}=\{x\in R^{1}\,|\,L^{x}_{t}=\sup_{y}L^{y}_{t}\},\label{rk.30}
\end{equation}
which we call the set of favorite points at time $t$. At any time $t$ there may be more than one favorite point. Let
\begin{equation}
V_{t}=\inf\{|x|\,|\, x\in \mathcal{V}_{t}\}.\label{rk.31}
\end{equation}
$V_{t}$ is a stochastic process in $t$. Does $\lim_{t\to\ff}V_{t}=\ff$? If so, how fast does $V_{t}$ grow? The generalized second Ray-Knight theorem has been used to give information about the rate of growth of $V_{t}$ for the symmetric stable processes, see Bass, Eisenbaum and Shi, \cite{BES} and the notes at the end of this chapter. We will illustrate this for the case of Brownian motion, although for this case one can avoid use of  the Ray-Knight theorem. Furthermore, in order not to get bogged down in details we consider only the following result:
\begin{equation}
\liminf_{t\to\ff}{\log^{\ga} t \over \sqrt{t}}V_{t}=\ff,\hspace{.2 in} P^{0} \hspace{.2 in}a.s.\label{rk.32}
\end{equation}
for any $\ga>6$. Note that the law of the iterated logarithm says that 
\[\lim\sup_{t\rar\ff}{B_{t} \over \sqrt{2t\log\log t}}=1, \hspace{.2 in} P^{0} \hspace{.2 in}a.s.\]  
so that in some sense the favorite points are near the boundary of the Brownian motion, $\approx \sqrt{t}$. 
Our techniques actually allow us to conclude that this holds for any $\ga>3$. It has been conjectured   that $\ga=1$ is the critical value. Furthermore, since our goal is only to illustrate how one can apply  the generalized second Ray-Knight theorem, we only discuss the proof for one direction of (\ref{rk.32}).

We use the generalized second Ray-Knight theorem to  prove the following. Let $h(t)=t/\log^{5} t$. Then  
\begin{equation}
\lim_{t\to\ff}\sup_{\{x\,|\,|x|\leq h(t)\}}{\log^{2} t \(L^{x}_{\tau(t)}-t\)\over t}=0,\hspace{.2 in} P^{0} \hspace{.2 in}a.s.\label{rk.33}
\end{equation}
For this, we first fix $\la$ large and   bound
\begin{equation}
P^{0}\( \sup_{\{x\,|\,|x|\leq h(t)\}} L^{x}_{\tau(t)}-t\geq 2\la \).\label{rk.34}
\end{equation}
This is certainly bounded by the following, which  allows us the opportunity to use the generalized second Ray-Knight theorem in the second line:
\bea
&&
\leq P^{0}P_{\eta}\( \sup_{\{x\,|\,|x|\leq h(t)\}} \(L^{x}_{\tau(t)}+\eta^{2}_{x}/2-t\)\geq 2\la \)
\label{rk.35}\\
&&=P_{\eta}\( \sup_{\{x\,|\,|x|\leq h(t)\}} \((\eta_{x}+\sqrt{2t})^{2}/2-t\)\geq 2\la \)\nn\\
&&=P_{\eta}\( \sup_{\{x\,|\,|x|\leq h(t)\}} \(\eta^{2}_{x}/2+\sqrt{2t}\eta_{x}\)\geq 2\la \)\nn\\
&&\leq P_{\eta}\( \sup_{\{x\,|\,|x|\leq h(t)\}}  \eta^{2}_{x} \geq 2\la \)+P_{\eta}\( \sup_{\{x\,|\,|x|\leq h(t)\}} \eta_{x}\geq \la/\sqrt{2t} \).\nn
\eea
Using the fact that $\eta_{x}$ is just $\sqrt{2}$ times two-sided Brownian motion together with the reflection principle allows us to bound the last line by
\begin{equation}
Ce^{-\la/h(t)}+Ce^{-\la^{2}/4t \,h(t)},\label{rk.36}
\end{equation}
Taking $\la=\ep t/\log^{2} t$ we have shown that  
\begin{equation}
P^{0}\( \sup_{\{x\,|\,|x|\leq h(t)\}} {\log^{2} t \(L^{x}_{\tau(t)}-t\)\over t}\geq 2\ep \)\leq C/t^{\ep^{2}/4}.\label{rk.37}
\end{equation}
Using Borel-Cantelli on the sequence $t_{n}=n^{8/\ep^{2}}$ and then interpolating we find that the left hand side of (\ref{rk.33}) is less than $2\ep$ for any $\ep>0$, which establishes (\ref{rk.33}).

Now assume that we can show that for $t\in [\tau(r)^{-},\tau(r)]$ we have the following lower bound on the absolute maximum local time
\begin{equation}
\sup_{x}L_{t}^{x}>r+r/\log^{2} r.\label{rk.38}
\end{equation}
By (\ref{rk.33}), using the fact that for $t\in [\tau(r)^{-},\tau(r)]$, $L_{t}^{0}\leq L_{\tau(r)}^{0}=r$
we have that for large $t$
\begin{equation}
\sup_{\{x\,|\,|x|\leq h(L_{t}^{0})\}} L_{t}^{x}\leq   \sup_{\{x\,|\,|x|\leq h(r)\}}L^{x}_{\tau(r)}\leq r+r/\log^{2} r.\label{rk.39}
\end{equation}
Comparison of this with (\ref{rk.38}) shows that 
\begin{equation}
V_{t}\geq h(L_{t}^{0}).\label{rk.40}
\end{equation}
We show below that for any $\ep>0$
\begin{equation}
\lim_{t\to\ff}{(\log t)^{1+\ep}  L^{0}_{t} \over \sqrt{t}}=\ff,\hspace{.2 in} P^{0} \hspace{.2 in}a.s.\label{rk.41}
\end{equation}
and together with (\ref{rk.40}) this gives  the lower bound in (\ref{rk.32}).

Before proving  (\ref{rk.41}), we observe that (\ref{rk.38}) can be obtained from 
\begin{equation}
\lim_{t\to\ff}\sup_{x}{\log^{2} t \(L^{x}_{\tau(t)}-t\)\over t}=\ff,\hspace{.2 in} P^{0} \hspace{.2 in}a.s.\label{rk.42}
\end{equation}
which can be shown by another application of the generalized second Ray-Knight theorem. However, for these lecture notes, one illustration is enough!

In order to prove (\ref{rk.41}) we need some basic facts about the inverse local time $\tau(t)$.
Since
\begin{equation}
\tau(t+s)=\tau(s)+\tau(t)\circ\th_{\tau(s)},\label{ta.1}
\end{equation}
it follows using the strong Markov property and the fact that $X_{\tau(s)}=0$ that
\bea
f(t+s)=:P^{0}\(e^{-\bb\tau(t+s)}\)&=&P^{0}\(e^{-\bb \tau(s)}\(e^{-\bb \tau(t)}\)\circ\th_{\tau(s)}\)\label{ta.2}\\
&=&P^{0}\(e^{-\bb \tau(s)}\)P^{0}\(e^{-\bb \tau(t)}\)=f(s)f(t).\nn
\eea
Since for $\bb>0$, $f(t)$ is decreasing, bounded by $1$ and right continuous, we must have
$f(t)=e^{-t\,v(\bb)}$ for some function $v(\la)$ which we now evaluate.

Note first that for any function $g$
\begin{equation}
\int_{0}^{\ff}g(t)\,dL_{t}^{0}=\int_{0}^{\ff}g(\tau(s))\,ds.\label{ta.3}
\end{equation}
To see this it suffices to verify it for functions of the form $g(t)=1_{\{[0,r]\}}(t)$, for which (\ref{ta.3})
is the statement that $L_{r}^{0}=|\{s\,|\,\tau(s) \leq r\}|$ which is easily verified. Then, using (\ref{ta.3})
\begin{eqnarray}
{1 \over v(\bb)}=\int_{0}^{\ff}f(s)\,ds&=&P^{0}\(\int_{0}^{\ff}e^{-\bb \tau(s)}\,ds\)
\label{ta.4}\\
&=& P^{0}\(\int_{0}^{\ff}e^{-\bb t}\,dL_{t}^{0}\) =u^{\bb}(0,0) \nonumber
\end{eqnarray}
by (\ref{rk.19}). Thus we have
\begin{equation}
P^{0}\(e^{-\bb \tau(t)}\)=e^{-t/u^{\bb}(0,0)}.\label{ta.5}
\end{equation}
In particular for Brownian motion, by (\ref{rk.14})
\begin{equation}
P^{0}\(e^{-\bb \tau(t)}\)=e^{-t\sqrt{2\bb}}.\label{ta.6}
\end{equation}

We next use this to show that  for Brownian motion
\begin{equation}
\limsup_{r\to\ff} {  \tau(r) \over r^{2}\,\log^{2+\ep}r}=0,\hspace{.2 in} P^{0} \hspace{.2 in}a.s.\label{ta.7}
\end{equation}
from which we will then derive (\ref{rk.41}). To prove (\ref{ta.7}) note that
\begin{eqnarray}
P^{0}\(\tau(r)\geq x\)&=&P^{0}\(1-e^{-  \tau(r)/x}\geq 1-e^{-1}\)
\label{ta.8}\\
&\leq &  {1 \over  1-e^{-1}} P^{0}\(1-e^{-  \tau(r)/x} \)\nn\\
&=&c\(1-e^{-r\sqrt{2/x}}\)\leq cr/\sqrt{x}.\nonumber
\end{eqnarray}
Thus
\begin{equation}
P^{0}\(\tau(r)\geq r^{2}\,\log^{2+\ep}r\)\leq c/\log^{1+\ep/2},\label{ta.9}
\end{equation}
so that taking $r_{n}=e^{n}$ by Borel-Cantelli
\begin{equation}
\limsup_{n\to\ff} {  \tau(r_{n}) \over r_{n}^{2}\,\log^{2+\ep}r_{n}}=0,\hspace{.2 in} P^{0} \hspace{.2 in}a.s.\label{ta.10}
\end{equation}
and (\ref{ta.7}) follows by interpolation.

Finally, (\ref{ta.7}) says that it takes less than $r^{2}\,\log^{2+\ep}r$ for the local time at $0$
to reach the level $r$, so that for large $r$
\begin{equation}
L^{0}_{r^{2}\,\log^{2+\ep}r}\geq r.\label{ta.11}
\end{equation}
Taking $r=\({t \over \log^{2+\ep}t}\)^{1/2}$ leads to (\ref{rk.41}).\qed

\subsection{Proof of the moment formula for $L^{y}_{\tau(\la)}$}

We first prove that
\begin{equation}
E^{x}_{\la}\(L^{y}_{\tau(\la)}\)=u_{T_{0}}(x,y)+\al.\label{rk.60}
\end{equation}

Since $\tau(t)$ cannot grow until the process first reaches $0$ we have $\tau(\la)=T_{0}+\tau(\la)\circ \th_{T_{0}}$. Hence for any $\bb>0$,
\begin{eqnarray}
 E^x_{
\la}\(\int_{\tau(\la)}^\ff e^{-\bb t}\,dL^{y}_t\) &=&E^x_{
\la}\(\int_{T_0+\tau(\la)\circ\th_{T_0}}^\ff e^{-\bb t}\,dL^{y}_t\)
\label{mp4.29z} \\ &=& E^x_{ \la}\(e^{ -\bb T_0}\int_{\tau(\la)\circ\th_{T_0}}^\ff e^{-\bb t}\,dL^{y}_t\)
\nonumber \\&=&E^x_{ \la}\(e^{ -\bb T_0}\(\int_{\tau(\la)}^\ff e^{-\bb t}\,dL^{y}_t\)\circ\th_{T_0}\)
\nonumber\\ &=&E^x\(e^{ -\bb T_0}E^{ X_{
T_0}}_{
\la}\(\int_{\tau(\la)}^\ff e^{-\bb t}\,dL^{y}_t\)\)
\nonumber\\ &=&E^x\(e^{ -\bb T_0}\)E^{ 0}_{
\la}\(\int_{\tau(\la)}^\ff e^{-\bb t}\,dL^{y}_t\)
\nonumber
\end{eqnarray}
using the strong Markov property at $T_{0}$, the fact that $P^{x}(T_0<\ff)=1 $ and that  $X_{ T_0}=0$. Similarly,
\begin{eqnarray}
\lefteqn{   E^{ 0}_{ \la}\(\int_{\tau(\la)}^\ff e^{-\bb t}\,dL^{y}_t\)
\label{mp4.29y} }\\ &&\qquad  =E^0_{ \la}\(1_{ \{
\tau(\la)<\ff\}}\int_{\tau(\la)}^\ff e^{-\bb t}\,dL^{y}_t\)
\nonumber \\ &&\qquad  =E^0_{ \la}\(e^{ -\bb \tau(\la)}1_{ \{
\tau(\la)<\ff\}}\(\int_{0}^\ff e^{-\bb t}\,dL^{y}_t\)\circ\th_{\tau(\la)}\)
\nonumber\\ &&\qquad  =E^0_{ \la}\(e^{ -\bb \tau(\la)}1_{ \{
\tau(\la)<\ff\}}E^{ X_{ \tau(\la)}}\(\int_{0}^\ff e^{-\bb t}\,dL^{y}_t\)\)
\nonumber\\ &&\qquad  =E^0_{ \la}\(e^{ -\bb \tau(\la)}\)E^{
0}\(\int_{0}^\ff e^{-\bb t}\,dL^{y}_t\)
\nonumber
\end{eqnarray} using the strong Markov property at $ \tau(\la)$, and  the fact that     $X_{ \tau(\la)}=0$ on $\tau(\la)<\ff$.
     Combining (\ref{mp4.29z}) and (\ref{mp4.29y}), we see that
\be E^x_{
\la}\(\int_{\tau(\la)}^\ff e^{-\bb
t}\,dL^{y}_t\)=E^x\(e^{ -\bb T_0}\)E^0_{ \la}\(e^{ -\bb
\tau(\la)}\)u^\bb(0,y).
\ee

     Using this and proceeding exactly as in (\ref{rk.20})--(\ref{rk.23}),
we see that
\begin{eqnarray} u_{\tau(\la)}(x,y)  &=&\lim_{\bb\rar 0}\lc
u^\bb(x,y)-{u^\bb(x,0)u^\bb(0,y)\over u^\bb(0,0) }\rc
\label{mp4.87}\\ && +\lim_{\bb\rar 0}{u^\bb(x,0)u^\bb(0,y)\over
u^\bb(0,0) }(1-E_\la^0(e^{-\bb
\tau(\la)})).\nn
\end{eqnarray} 
 By (\ref{rk.13}),
\begin{equation}
\lim_{\bb\rar 0}\lc u^\bb(x,y)-{u^\bb(x,0)u^\bb(0,y)\over u^\bb(0,0)
}\rc=u_{\mbox{\tiny $T_0$}}(x,y).\label{mp4.88}
\end{equation}
      Also, by (\ref{rk.21}),
\bea &&
\lim_{\bb\rar 0}{u^\bb(x,0)u^\bb(0,y)\over u^\bb(0,0)
}(1-E_\la^0(e^{-\bb
\tau(\la)}))\label{mp4.28}\\ &&\qquad=\lim_{\bb\rar 0}E^x(e^{-\bb
T_0})E^y(e^{-\bb T_0})u^\bb(0,0)(1-E_\la^0(e^{-\bb
\tau(\la)}))
\nn\\ &&\qquad=P^x(T_0<\ff)P^y(T_0<\ff)\lim_{\bb\rar
0}u^\bb(0,0)(1-E_\la^0(e^{-\bb
\tau(\la)}))\nn\\ &&\qquad=\lim_{\bb\rar
0}u^\bb(0,0)(1-E_\la^0(e^{-\bb
\tau(\la)})),\nn
\eea and, by (\ref{ta.5}),
\begin{eqnarray}
      && u^\bb(0,0)(1-E_\la^0(e^{-\bb
\tau(\la)})) \label{mp4.28a}\\
      &&\qquad =u^\bb(0,0)(1-E_{ \la}(e^{-\la/u^{ \bb}( 0,0)}))   \nn\\
      &&\qquad =u^\bb(0,0)(1-{ 1/\al\over 1/\al+1/u^{ \bb}( 0,0)})={1\over
1/\al+1/u^{
\bb}( 0,0)}. \nn
\end{eqnarray}
Since $\lim_{\bb\rar 0}u^{ \bb}( 0,0)=\ff$ we get (\ref{rk.60}).

Combing    (\ref{rk.60})  with the definition (\ref{rk.6}) we have 
\bea  u_{\tau (\la)}(x,y):=u_{T_{0}}(x,y)+\al&=&E^x_\la\(L^y_{\tau (\la)}\)\label{bm8.8d}\\ 
 &=&E^x_\la\(
\int_0^{\ff}1_{\{\tau (\la)> t\}}\,dL^{y}_{t}\)\nn\\&=&
E^x_\la\(
\int_0^{\ff}1_{\{\la> L_{t}^{0}\}}\,dL^{y}_{t}\)\nn\\ &=& E^x \(
\int_0^{\ff}P_{ \la}(\la> L_{t}^{0})\,dL^{y}_{t}\)\nn\\ &=&E^x\(
\int_0^\ff e^{- L_{t}^{0}/\al}\,dL^{y}_{t}\).\nn
\eea

We now prove (\ref{rk.6}). 
 We have
\begin{eqnarray} E_{\la}^x\(\prod_{i=1}^n
L^{y_i}_{\tau (\la)} \)&=&E_{\la}^x\(\prod_{i=1}^n
\int_0^{\tau (\la)}
\,dL^{y_i}_{t_i}\)\label{bm2.18}\\ &=&   \sum_{\pi\in \PP_{n}} E_{\la}^x\(
\int_{\{0<t_1<\ldots<t_n<\tau (\la)\}}\prod_{i=1}^n\,dL^{y_{\pi_i}}_{t_i}\).\nn
\end{eqnarray} 

Let $y_{\pi_i}=z_i$, $i=1,\ldots,n$. Note that
\be
\int_{\{0<t_1<\ldots<t_n<\tau (\la)\}}\prod_{i=1}^n\,dL^{z_{i}}_{t_i}
=\int_0^{\tau (\la)}\!\!\int_{t_1}^{\tau (\la)}\cdots
\int_{t_{n-1}}^{\tau (\la)} \prod_{i=1}^n\,dL^{z_{i}}_{t_i}.\label{bm2.18qa}
\ee    Therefore, setting
\begin{equation}
F(t,\tau (\la))=\int_{\{t<t_2<\ldots<t_n<\tau (\la)\}}\prod_{i=2}^n\,dL^{z_{i}}_{t_i},
\label{bm2.22p}
\end{equation} we have
\bea
&&
E^x_{\la}\(\int_{\{0<t_1<\ldots<t_n<\tau (\la)\}}\prod_{i=1}^n\,dL^{z_{i}}_{t_i}\)
=E^x_{\la}\(
\int_0^{\tau (\la)} F(t,\tau (\la))\,dL^{z_1}_{t}\)\nn\\
&&\hspace{1.5 in}=E^x_{\la}\(
\int_0^{\ff} 1_{\{\la>L_{t}^{0}\}}F(t,\tau (\la))\,dL^{z_1}_{t}\).\label{bm2.21}
\eea

Since  $\la>L_{t}^{0}$ implies that
$\tau(\la)=t+\tau(\la-L_{t}^{0})\circ\th_t $,
\bea\lefteqn{ E_\la\(1_{\{\la> L_{t}^{0}\}} F(t,\tau (\la))\)} \\ &&=
E_\la\(1_{\{\la> L_{t}^{0}\}} F(t,t+\tau (\la-L_{t}^{0})\circ\th_t)\)\nn\\
&&= \int_{L_{t}^{0}}^\ff F(t,t+\tau (y-L_{t}^{0})\circ\th_t)e^{-  y/\al}\,dy/\al\nn\\
&&=e^{-  L_{t}^{0}/\al}\int_{0}^\ff F(t,t+\tau (y )\circ\th_t)e^{- 
y/\al}\,dy/\al\nn\\
&&=e^{-  L_{t}^{0}/\al}E_\la\( F(t,t+\tau (\la )\circ\th_t)\).\nn
\eea

Using this in (\ref{bm2.21})
\bea && E^x_\la\( \int_0^{\ff}1_{\{\la>
L_{t}^{0}\}}F(t,\tau (\la))\,dL^{z_1}_{t}\)\label{bm8.8a}\\ &&\qquad
=E^x_\la\(
\int_0^{\ff}e^{-  L_{t}^{0}/\al}
F(t,t+\tau (\la)\circ\th_t)\,dL^{z_1}_{t}\)\nn\\ &&\qquad =E^x_\la\(
\int_0^{\ff}e^{-  L_{t}^{0}/\al} F(0,\tau (\la))\circ\th_t\,dL^{z_1}_{t}\),\nn
\eea 
where the last equality uses the additivity of local times.

Let $\tau_{z_{1}}(s)$ denote the right continuous local time for $z_{1}$. Using the analogue of (\ref{ta.3}) for $\tau_{z_{1}}(s)$  and then the strong Markov property at $\tau_{z_{1}}(s)$ we have
\bea && E^x_\la\(
\int_0^{\ff}e^{-  L_{t}^{0}/\al}
F(0,\tau (\la))\circ\th_t\,dL^{z_1}_{t}\)\label{bm8.8ff}\\ &&= E^x_\la\(
\int_0^{\ff}e^{-  L_{\tau_{z_{1}}(s)}^{0}/\al}
F(0,\tau (\la))\circ\th_{\tau_{z_{1}}(s)}\,ds \)\nn\\ &&= \int_0^{\ff}E^x_\la\(
e^{-  L_{\tau_{z_{1}}(s)}^{0}/\al}
F(0,\tau (\la))\circ\th_{\tau_{z_{1}}(s)}\)\,ds \nn\\ &&= \int_0^{\ff}E^x\(
e^{-  L_{\tau_{z_{1}}(s)}^{0}/\al}E_\la^{z_{1}}\(
F(0,\tau (\la))\)\)\,ds \nn\\&&\qquad
    =E_\la^{z_1}\(F(0,\tau (\la))\) E^x\(
\int_0^\ff e^{-  L_{t}^{0}/\al}\,dL^{z_1}_{t}\)\nn\\ &&\qquad
    = E_\la^{z_1}\(F(0,\tau (\la))\)u_{ \tau (\la)}(x,z_{1}),\nn
\eea
    where, for the next to  last equation, we use  (\ref{ta.3}) and for the last equation, we use  (\ref{bm8.8d}).

    Combining (\ref{bm2.18})--(\ref{bm8.8ff}) we see that

\bea  && E^x_\la\(
\int_{\{0<t_1<\ldots<t_n<\tau (\la)\}}\prod_{i=1}^n\,dL^{z_{i}}_{t_i}\)
\nn\\ &&\qquad=u_{ \tau (\la)}(x,z_1)
E^{z_1}_\la\(\int_{\{0<t_2<\ldots<t_n<\tau (\la)\}}
\prod_{i=2}^n\,dL^{z_{i}}_{t_i}\).\nn
\eea Iterating this argument we get
\bea  && E^x_\la\(
\int_{\{0<t_1<\ldots<t_n<\tau (\la)\}}\prod_{i=1}^n\,dL^{z_{i}}_{t_i}\)
\label{bm2.20qq}\\ &&\qquad =u_{ \tau (\la)}(x,z_{1})u_{
\tau (\la)}(z_{1},z_{2})\cdots u_{ \tau (\la)}(z_{n-1},z_{n}).\nn
\eea  
 Summing over
$\pi$ we get  (\ref{rk.6}).

{\it     The proof of  the  generalized second Ray-Knight theorem given here is similar to that   in \cite[Section 8.3.3]{book}.   Chapter 11 gives a full treatment of favorite points for Brownian motion and stable processes.

We have assumed that our Markov process is recurrent. For the transient case see \cite[Theorem 8.2.3]{book}.}

\section{Loop soups}

Until now we have considered only symmetric Markov processes. This was natural since our Isomorphism theorems related Markov local times to the squares of an associated Gaussian process whose covariance was the potential density of our Markov process, and covariance functions are always symmetric. If we want to have an Isomorphism theorem for the local times of a not necessarily symmetric Markov processes we will need to find a substitute for Gaussian squares. The route we take is long, but quite interesting. In the end it should help remove some of the mystery of Isomorphism theorems, even in the symmetric case. The mystery I refer to  is this: Even after all the proofs we have seen, why, intuitively, should 
Gaussian squares be related to Markov local times?

\subsection{The loop measure}

Once agin, we assume a Markov process $X\in S$ with transition densities $p_{t}(x,y)$ with respect to a measure $m$ on S. But we do not assume that $p_{t}(x,y)$ is symmetric. Our first step is to define bridge measures $ Q_{t}^{x,y}$ for $X$. Consider
 \begin{equation}
M_{s}=p_{t-s}(X_{s},x).\label{pe.1}
\end{equation}
Let us show that $M_{s}$ is a martingale, $0\leq s\leq t$. Note that
 \begin{equation}M_{s}=p_{t-s}(X_{s},x)= p_{t-s}(X_{s-r},x)\circ\theta_{r}.
\label{pe.2}
\end{equation}
Hence, using the Markov property 
 \bea
 &&E^{x}\(M_{s}\,|\,\mathcal{F}_{r}\)=E^{X_{r}}\(p_{t-s}(X_{s-r},x) \)\label{pe.3}\\
 &&
=\int p_{s-r}(X_{r},z) p_{t-s}(z,x)\,dm(z)=M_{r}.\nn
 \eea
It follows that if we set 
 \begin{equation}
 Q_{t}^{x,y}\(G\)=P^{x}\( G  \,M_{s} \)=P^{x}\( G  \,p_{t-s}(X_{s},y) \), 
 \label{pe.4}
\end{equation}
 for all $G\in \mathcal{F}_{s}$ with $s<t$,
 then $Q_{t}^{x,y}$ is well-defined. That is, if in fact $G\in \mathcal{F}_{r}$ with $r<s<t$, then 
 we obtain the same value for $ Q_{t}^{x,y}\(G\)$ if we use
  \begin{equation}
 Q_{t}^{x,y}\(G\)=P^{x}\( G  \,M_{r} \), 
 \label{pe.5}
\end{equation}
which follows from the fact that $M_{s}$ is a martingale, $0\leq s\leq t$.
We note that $Q_{t}^{x,y}$ extends naturally  to $\mathcal{F}^{-}_{t}$.

 Given bridge measures we can now define the loop measure
 \begin{equation}
 \mu (F)=\int_{0}^{\ff}{1\over t}\int\,Q_{t}^{x,x}\(F\circ  k_{t} \)\,dm (x)\,dt
 \label{pe.6}
\end{equation}
 for all $ F\in \mathcal{F}$, 
      where  
 $k_{t}$ is the killing operator:  $k_{t}\om(s)=\om(s)$   if $s<t$,   and  $ k_{t}\om(s)=\De$   if $s\geq t$.   
 
We need to see how to compute with $\mu$. Our goal is to obtain the following  moment formula for local times 
under $\mu$. 
\bea
\mu\(\prod_{j=1}^{k}L_{\ff}^{x_{j}}\)&=&\mbox{cy}([1,k]) \label{pe.loop}\\&=& \sum_{\pi\in \mathcal{P}^{\odot}_{k}}  u(x_{\pi(1)},x_{\pi(2)})\cdots   u(x_{\pi(k-1)},x_{\pi(k)}) u(x_{\pi(k)},x_{\pi(1)}).\nn
\eea

Proof: 
If $0<t_{1}<t_{2}<\cdots<t_{k}$ it follows from the definition of the killing operator 
 $k_{t}$ that
 \begin{equation}
  \prod_{j=1}^{k}f_{j}(X_{t_{j}})\circ  k_{t}=1_{\{t>t_{k}\}} \prod_{j=1}^{k}f_{j}(X_{t_{j}}).
  \label{pe.7}
\end{equation}
Hence from the definition of the  bridge measure $Q_{t}^{x,x}$
\begin{eqnarray}
&&\hspace{-.3 in}Q_{t}^{x,x}\(\prod_{j=1}^{k}f_{j}(X_{t_{j}})\circ  k_{t} \)  =1_{\{t>t_{k}\}}  \int   p_{t_{1}}(x, y_{1})   f_{1}(y_{1}) p_{t_{2}-t_{1}}(y_{1},y_{2})f_{2}(y_{2})\cdots \nonumber\\
&&\hspace{1 in}\cdots  p_{t_{k}-t_{k-1}}(y_{k-1},y_{k})f_{k}(y_{k})p_{t-t_{k}}(y_{k},x)\,dm(y_{1})
\cdots \,dm (y_{k})\nn
\end{eqnarray}
 Integrating with respect to $dm (x)$ and using $\int p_{t-t_{k}}(y_{k},x)  p_{t_{1}}(x, y_{1})\,dm (x)=p_{t_{1}+t-t_{k}}(y_{k},y_{1})$
we obtain\begin{eqnarray}
&&\hspace{-.3 in}
\mu\(   \prod_{j=1}^{k}f_{j}(X_{t_{j}})\) =\int_{t_{k}}^{\ff}{1\over t} \int f_{1}(y_{1}) p_{t_{2}-t_{1}}(y_{1},y_{2})f_{2}(y_{2})\cdots   \label{pe.8}\\
&&\cdots  p_{t_{k}-t_{k-1}}(y_{k-1},y_{k})f_{k}(y_{k})p_{t_{1}+t-t_{k}}(y_{k},y_{1})\,dm(y_{1})\cdots \,dm (y_{k}).\nn
\end{eqnarray}
 This does not look very enlightening, but we plough ahead and  integrate over time to obtain
\begin{eqnarray}
 && \hspace{-.3 in}
\mu\( \int_{\{0\leq t_{1}\leq \cdots \leq t_{k-1}\leq t_{k} <\ff\}}  \prod_{j=1}^{k}f_{j}(X_{t_{j}})\,dt_{j}\)   \label{pe.9}\\
&& \hspace{-.3 in}=\int_{\{0\leq t_{1}\leq \cdots \leq t_{k} \leq t<\ff\}}  {1 \over t} \int f_{1}(y_{1}) p_{t_{2}-t_{1}}(y_{1},y_{2})f_{2}(y_{2})\cdots \nonumber\\
&&\hspace{.1 in}\cdots  p_{t_{k}-t_{k-1}}(y_{k-1},y_{k})f_{k}(y_{k})p_{t_{1}+t-t_{k}}(y_{k},y_{1}) \prod_{j=1}^{k} \,dm(y_{j})\,dt_{j}\,dt.\nn
\end{eqnarray}
Note that we are integrating over $k+1$ time variables, $t_{1}\leq \cdots \leq t_{k} \leq t$. We make the following change of variables: $r_{1}=t_{1}+t-t_{k}, r_{2}=t_{2}-t_{1},\ldots r_{k}=t_{k}-t_{k-1}$, and retain $t_{1}$ as our $k+1$'st time variable. The range of integration is $[0.\ff]$ for all $r_{j}$
 and $t_{1}\leq r_{1}$. We also note that $r_{1}+ r_{2}+\cdots +r_{k}=t$ so that (\ref{pe.9}) becomes
\begin{eqnarray}
&& \hspace{-.3 in}
\mu\( \int_{\{0\leq t_{1}\leq \cdots \leq t_{k-1}\leq t_{k}<\ff\}}  \prod_{j=1}^{k}f_{j}(X_{t_{j}})\,dt_{j}\) \label{pe.10}\\
&& \hspace{-.3 in}=\int_{R^{k}_{+}}  {1\over r_{1}+\cdots +r_{k}}\( \int f_{1}(y_{1}) p_{r_{2}}(y_{1},y_{2})f_{2}(y_{2})\cdots \right.\nonumber\\
&&\hspace{.3 in}\left.\cdots  p_{r_{k}}(y_{k-1},y_{k})f_{k}(y_{k})p_{r_{1}}(y_{k},y_{1}) \prod_{j=1}^{k} \,dm(y_{j})\)\(\int_{0}^{r_{1}} 1\,dt_{1}\)\prod_{j=1}^{k} \,dr_{j}\nn
\\
&&\hspace{-.3 in}=\int_{R^{k}_{+}}     {r_{1}\over r_{1}+\cdots +r_{k}} \int f_{1}(y_{1}) p_{r_{2}}(y_{1},y_{2})f_{2}(y_{2})\cdots \nonumber\\
&&\hspace{1 in}\cdots  p_{r_{k}}(y_{k-1},y_{k})f_{k}(y_{k})p_{r_{1}}(y_{k},y_{1})  \prod_{j=1}^{k} \,dm(y_{j})\,dr_{j}.\nn
\end{eqnarray}
 
Using the fact that $R^{k}_{+}=\cup_{\pi\in \mathcal{P}_{k}}\{0\leq t_{\pi(1)}\leq \cdots \leq t_{\pi(k)}\leq t_{\pi(k-1)}<\ff\}$
we then obtain\begin{eqnarray}
&& \hspace{-.3 in}
\mu\( \prod_{j=1}^{k}\( \int_{0}^{\ff}  f_{j}(X_{t})\,dt\)\)\label{pe.11}\\
&&
 \hspace{-.3 in}=\sum_{\pi\in \mathcal{P}_{k}}\mu\( \int_{\{0\leq t_{1}\leq \cdots \leq t_{k-1}\leq t_{k}<\ff\}}  \prod_{j=1}^{k}f_{\pi(j)}(X_{t_{j}})\,dt_{j}\) \nn\\
&& \hspace{-.3 in}=\int    {r_{1} \over r_{1}+\cdots +r_{k}}\sum_{\pi\in \mathcal{P}_{k}} \int f_{\pi(1)}(y_{1}) p_{r_{2}}(y_{1},y_{2})f_{\pi(2)}(y_{2})\cdots \nonumber\\
&&\cdots  p_{r_{k}}(y_{k-1},y_{k})f_{\pi(k)}(y_{k})p_{r_{1}}(y_{k},y_{1})  \prod_{j=1}^{k} \,dm(y_{j})\,dr_{j}\nn\\
&& \hspace{-.3 in}=\int    {r_{1} \over r_{1}+\cdots +r_{k}}h(r_{1},r_{2},\ldots,r_{k} )\prod_{j=1}^{k} \,dr_{j}, \nn 
\end{eqnarray}
where 
\begin{eqnarray}
\lefteqn{
h(r_{1},r_{2},\ldots,r_{k} )\label{pe.12}}\\
&& \hspace{-.4 in}=\sum_{\pi\in \mathcal{P}_{k}} \int  p_{r_{2}}(y_{1},y_{2})\cdots  p_{r_{1}}(y_{k},y_{1})  \prod_{j=1}^{k}f_{\pi(k)}(y_{k}) \,dy.\nn 
\end{eqnarray}
 The basic idea we now use  is that since $h$ involves a sum over permutations, there is no longer anything special about $r_{1}$, which will allow us to eliminate the factor $   {r_{1} \over r_{1}+\cdots +r_{k}}$ and end up with a nice formula. In more detail, observe that
\begin{equation}
h(r_{1},r_{2},\ldots,r_{k} )=h(r_{2},r_{3},\ldots,r_{1} ).\label{pe.13}
\end{equation}
Hence, first changing variables in (\ref{pe.11}), and then using (\ref{pe.13}) we obtain
\begin{eqnarray} \hspace{-.3 in}
\mu\( \prod_{j=1}^{k}\( \int_{0}^{\ff}  f_{j}(X_{t})\,dt\)\)&=&\int    {r_{2} \over r_{1}+\cdots +r_{k}}h(r_{2},r_{3},\ldots,r_{1} )\prod_{j=1}^{k} \,dr_{j}\nn\\
&=& \int    {r_{2} \over r_{1}+\cdots +r_{k}}h(r_{1},r_{2},\ldots,r_{k} )\prod_{j=1}^{k} \,dr_{j}.\label{pe.11a} 
\end{eqnarray}
Adding together similar expressions where $r_{2}$ in the numerator  is replaced in turn by  $r_{3},\cdots ,r_{k}$ we have shown that 
\bea 
&&\hspace{-.35 in}
\mu\( \prod_{j=1}^{k}\( \int_{0}^{\ff}  f_{j}(X_{t})\,dt\)\)={1 \over k}\int   h(r_{1},r_{2},\ldots,r_{k} )\prod_{j=1}^{k} \,dr_{j}\label{pe.13s}\\
&&\hspace{-.35 in} ={1 \over k} \sum_{\pi\in \mathcal{P}_{k}} \int f_{\pi(1)}(y_{1}) u(y_{1},y_{2})f_{\pi(2)}(y_{2})\cdots   u(y_{k-1},y_{k})f_{\pi(k)}(y_{k})u(y_{k},y_{1})  \prod_{j=1}^{k} \,dm(y_{j}),\nonumber
\eea
Letting $f_{j}=f_{x_{j},\de}$ and taking the limit $\de\rar 0$ we obtain a simple moment formula:
\be
\mu\(\prod_{j=1}^{k}L_{\ff}^{x_{j}}\) ={1 \over k} \sum_{\pi\in \mathcal{P}_{k}}  u(x_{\pi(1)},x_{\pi(2)})\cdots   u(x_{\pi(k-1)},x_{\pi(k)}) u(x_{\pi(k)},x_{\pi(1)}). \label{pe.14}
\ee
The product of $u$'s on the right is invariant under the  $k$ rotations $( 1,2,\ldots,k)\to ( i,i+1,\ldots, i+k)$ mod $k$. Thus we have
\be
\mu\(\prod_{j=1}^{k}L_{\ff}^{x_{j}}\) = \sum_{\pi\in \mathcal{P}^{\odot}_{k}}  u(x_{\pi(1)},x_{\pi(2)})\cdots   u(x_{\pi(k-1)},x_{\pi(k)}) u(x_{\pi(k)},x_{\pi(1)}), \label{pe.14a}
\ee
which is (\ref{pe.loop}).\qed

Fixing a point $x_{0}$ then gives us  
\be
\hspace{-.1 in}\mu\(L_{\ff}^{x_{0} }\prod_{j=1}^{k}L_{\ff}^{x_{j}}\)  = \sum_{\pi\in \mathcal{P}_{k}} u(x_{0},x_{\pi(1)}) u(x_{\pi(1)},x_{\pi(2)})\cdots   u(x_{\pi(k-1)},x_{\pi(k)}) u(x_{\pi(k)},x_{0}),   \label{pe.15}
\ee
or, using (\ref{mp.14}),
\begin{equation}
\mu\(L_{\ff}^{x_{0} }\prod_{j=1}^{k}L_{\ff}^{x_{j}}\)=\mbox{cy}([0,k])=Q^{x_{0}, x_{0}}  \(\prod_{j=1}^{k}L_{\ff}^{x_{j}}\). \label{pe.16}
\end{equation}
Recalling the role that $Q^{x, x} $ played in Dynkin's isomorphism theorem in the symmetric case, we can feel we are getting closer to an isomorphism theorem in the non-symmetric case. We  need to recall some basic facts about Poisson  processes.

But first we make a slight improvement on (\ref{pe.16}). Since we have already seen that total local times are exponentially integrable under $Q^{x, x}$, and this does not depend on symmetry, (\ref{pe.16}) implies that
 \begin{equation}
\mu\(L_{\ff}^{x } F\(L_{\ff}^{x_{i}}\)  \)=Q^{x, x}  \(  F\(L_{\ff}^{x_{i}}\)   \). \label{pe.17}
\end{equation}

 \subsection{Poisson processes}

 Let $\mathcal{L}_{\al}$ be  a Poisson   process  on $\Om_{\De}$ with intensity measure $\al \mu$.  Thus, each realization of $\mathcal{L}_{\al}$ is a countable collection of points in $\Om_{\De}$, and if
 \begin{equation}
N(A):=\# \{\mathcal{L}_{\al}\cap A \},\hspace{.2 in}A\subseteq \Om_{\De} \label{pp.1}
\end{equation}
then
 \begin{equation}
 P_{\mathcal{L}_{\alpha}}\(N(A)=k\)={(\al \mu (A))^{k} \over k!}e^{-\al \mu (A)},
 \label{pp.2}
\end{equation}
and $N(A_{1}), \ldots,N(A_{k})$ are independent for disjoint $A_{1},\ldots, A_{k}$. 
For any bounded measurable functional $f$ on $\Om_{\De}$ let
\begin{equation}
N(f)= \sum_{\om\in \mathcal{L_{\al}}} f(\om),\label{pp.2p}
\end{equation}
so that $N(A)=N(1_{\{A\}})$.

 We will need three basic facts about our Poisson  process. The master formula for Poisson processes says that for any bounded measurable functional $f$ on $\Om_{\De}$
  \begin{equation}
E_{\mathcal{L}_{\alpha}}\(e^{N(f)}\)=\exp \(\al\(\int_{\Om_{\De}}\( e^{f(\om)}-1\)\,d\mu(\om)  \)\).
 \label{pp.3}
\end{equation}
Proof: A simple calculation using (\ref{pp.2}) shows that
\begin{equation}
E_{\mathcal{L}_{\alpha}}\(e^{z N(A)}\)=\sum_{k=0}^{\ff}e^{z k}{(\al \mu (A))^{k} \over k!}e^{-\al \mu (A)}=\exp\(\al \(e^{z}-1\)\mu (A)\).\label{pp.2a}
\end{equation}
If $A_{1}\cup\cdots\cup  A_{n}=\Om_{\De}$ is a partition of $\Om_{\De}$, and $f=\sum^{n}_{j=1}z_{j}1_{\{A_{j}\}}$ then
\begin{equation}
N(f)= \sum_{\om\in \mathcal{L_{\al}}} f(\om) =\sum^{n}_{j=1}z_{j}N(A_{j}),\label{pp.2b}
\end{equation}
so that by independence we have 
\bea
&&
E_{\mathcal{L}_{\alpha}}\(e^{N(f)}\)=\prod^{n}_{j=1} E_{\mathcal{L}_{\alpha}}\(e^{z_{j} N(A_{j})}\)     \label{pp.2c}\\
&&=\prod^{n}_{j=1}\exp\(\al \(e^{z_{j}}-1\)\mu (A_{j})\) = \exp \(\al\(\int_{\Om_{\De}}\( e^{f(\om)}-1\)\,d\mu(\om)  \)\),\nn
\eea
and (\ref{pp.3}) for general $f$ follows on taking limits.\qed

The second fact is the moment formula 
\begin{equation}
E_{\mathcal{L}_{\alpha}}\( \prod_{j=1}^{n}N(f_{j})  \)=\sum_{\cup_{i=1}^{\ell} B_{i}=[1,n]} \,\,\prod_{i=1}^{\ell}\,\al\,\mu\(\prod_{j\in B_{i}} f_{j}\).\label{pp.15}
\end{equation}

Proof: Since $N(\sum_{j=1}^{n}z_{j}f_{j})=\sum_{j=1}^{n}z_{j}N(f_{j})$, by the master formula
 \begin{equation}
E_{\mathcal{L}_{\alpha}}\(e^{ \sum_{j=1}^{n}z_{j}N(f_{j}) }\)=\exp \(\al\(\int_{\Om_{\De}}\( e^{\sum_{j=1}^{n}z_{j}f_{j}(\om)}-1\)\,d\mu(\om)  \)\). \label{pp.12}
\end{equation} 
Differentiating with respect to $z_{1}$  and then setting $z_{1}=0$ we obtain 
\bea
&&\hspace{-.4 in}
E_{\mathcal{L}_{\alpha}}\(N(f_{1})e^{ \sum_{j=2}^{n}z_{j}N(f_{j}) }\) \label{pp.13}\\
&&\hspace{-.4 in}=\al\(\int_{\Om_{\De}}  f_{1}e^{\sum_{j=2}^{n}z_{j}f_{j}(\om)}\,d\mu(\om)  \)\exp \(\al\(\int_{\Om_{\De}}\( e^{\sum_{j=2}^{n}z_{j}f_{j}(\om)}-1\)\,d\mu(\om)  \)\).\nn
\eea
Differentiate now with respect to $z_{2}$, using the product rule for the right hand side  and then setting $z_{2}=0$ we obtain 
\bea
&&\hspace{-.4 in}
E_{\mathcal{L}_{\alpha}}\(N(f_{1})N(f_{2})e^{ \sum_{j=3}^{n}z_{j} N(f_{j}) }\) \label{pp.14}\\
&&\hspace{-.4 in}=\al\(\int_{\Om_{\De}} f_{1}f_{2}e^{\sum_{j=3}^{n}z_{j}f_{j}(\om)}\,d\mu(\om)  \)\exp \(\al\(\int_{\Om_{\De}}\( e^{\sum_{j=3}^{n}z_{j}f_{j}(\om)}-1\)\,d\mu(\om)  \)\)\nn\\
&&\hspace{-.2 in}+\al\(\int_{\Om_{\De}}  f_{1}e^{\sum_{j=3}^{n}z_{j}f_{j}(\om)}\,d\mu(\om)  \)\al\(\int_{\Om_{\De}}f_{2}e^{\sum_{j=3}^{n}z_{j}f_{j}(\om)}\,d\mu(\om)  \)\nn\\
&&\hspace{2 in}\exp \(\al\(\int_{\Om_{\De}}\( e^{\sum_{j=3}^{n}z_{j}f_{j}(\om)}-1\)\,d\mu(\om)  \)\).\nn
\eea
By now it should be clear that iterating this leads to  (\ref{pp.15}).

Our last basic fact is the Palm formula which  says that for $f$ as above and $G$ a symmetric measurable function on $\Om_{\De}^{\ff}$
  \begin{equation}
E_{\mathcal{L}_{\alpha}} \(N(f)G(\mathcal{L}_{\alpha})\) =\al
\int E_{\mathcal{L}_{\alpha}} \(G(\om'\cup\,\mathcal{L}_{\alpha})\)f(\om')\,d\mu (\om'). \label{pp.4}
\end{equation}
Proof: By the master formula
  \begin{equation}
E_{\mathcal{L}_{\alpha}}\(e^{ zN(f)+N(g)} \)=\exp \(\al\(\int_{\Om_{\De}}\( e^{zf(\om)+g(\om)}-1\)\,d\mu(\om)  \)\).
 \label{pp.3d}
\end{equation}
Differentiation with respect to $z$ and then setting $z=0$ we obtain
  \bea
  &&
E_{\mathcal{L}_{\alpha}}\(N(f)\,\,e^{ N(g)} \) \label{pp.3d}\\
&&=\al\(\int_{\Om_{\De}}f(\om) e^{g(\om)}\,d\mu(\om)  \)\exp \(\al\(\int_{\Om_{\De}}\( e^{g(\om)}-1\)\,d\mu(\om)  \)\).
\nn
\eea
Using the master formula on the right most term we can write the right hand side as
  \bea
  &&
  \al\(\int_{\Om_{\De}}f(\om') e^{g(\om')}\,d\mu(\om')  \) E_{\mathcal{L}_{\alpha}}\(e^{ \sum_{\om\in \mathcal{L_{\al}}}g(\om) } \)  \label{pp.3e}\\
  &&=   \al\(\int_{\Om_{\De}}f(\om')  E_{\mathcal{L}_{\alpha}}\(e^{ \sum_{\om\in \om'\cup \mathcal{L_{\al}}}g(\om) } \)\,d\mu(\om')  \).
\nn
\eea
This proves (\ref{pp.4}) for the special case when $G(\om_{i})=e^{\sum_{i=1}^{\ff}g(\om_{i})}$, and this would actually be sufficient for our purposes, but in fact the general case follows from this.\qed

 \subsection{The  isomorphism theorem}
  We now use loop soups to prove a general isomorphism theorem for not necessarily symmetric Markov processes, which in the symmetric case,  with $\al=1/2$, is the Dynkin isomorphism theorem.
 Let  
 \begin{equation}
 \wh L_{\al}^{x}=N(L_{\ff}^{x})=\sum_{\om\in \mathcal{L_{\al}}} L_{\ff}^{x}(\om).
 \label{pp.5}
\end{equation}
 Using  (\ref{pp.15}) with $f_{j}=L_{\ff}^{x_{j}}$ and then  (\ref{pe.loop}) we obtain 
\bea
E_{\mathcal{L}_{\alpha}}\( \prod_{j=1}^{n}\wh L_{\al}^{x_{j}}  \)&=&\sum_{\cup_{i=1}^{\ell} B_{i}=[1,n]}\,\,\,\prod_{i=1}^{\ell} \al \,\mu\(\prod_{j\in B_{i}} L_{\ff}^{x_{j}}\)\label{pp.15r}\\
&=&\sum_{\cup_{i=1}^{\ell} B_{i}=[1,n]}\,\,\,\prod_{i=1}^{\ell}\al \,\mbox{cy}(B_{i}).\nn
\eea
We prove the following general isomorphism theorem
 \begin{equation}
E_{\mathcal{L}_{\al}} \(\wh L_{\al}^{x_{0}} \,F\(\wh L_{\al}^{x_{j}}\)  \)
=\al E_{\mathcal{L}_{\al}} Q^{x_{0},x_{0}}\(  F \( \wh L_{\al}^{x_{j}} +L_{\ff}^{x_{j}}\) \). \label{pp.9v}
\end{equation}
Proof: As before, it suffices to prove that
 \begin{equation}
E_{\mathcal{L}_{\al}} \(\wh L_{\al}^{x_{0}} \,\prod_{j=1}^{k}\wh L_{\al}^{x_{j}}  \)
=\al E_{\mathcal{L}_{\al}} Q^{x_{0},x_{0}}\(  \prod_{j=1}^{k} \( \wh L_{\al}^{x_{j}} +L_{\ff}^{x_{j}}\) \),\label{pp.9e}
\end{equation}
which we can write as 
  \begin{equation}
 E_{\mathcal{L}_{\al}} \(\prod_{j=0}^{k}\wh L_{\al}^{x_{j}}   \)=\al\sum_{A\subseteq [1,k]} Q^{x_{0},x_{0}}\(\prod_{i\in A}   L^{ x_{i}}_{\ff}   \)E_{\mathcal{L}_{\al}}\(\prod_{i\notin A}\wh L_{\al}^{x_{i}} \).\label{dit.3loop}
 \end{equation}
 Using (\ref{pp.15r}), our theorem is the claim that 
 \begin{eqnarray}
 &&\sum_{\cup_{i=0}^{\ell} B_{i}=[0,n]}\,\,\,\prod_{i=1}^{\ell}\al \,\mbox{cy}(B_{i})
 \label{}\\
 && =\al\sum_{A\subseteq [1,k]} Q^{x_{0},x_{0}}\(\prod_{i\in A}   L^{ x_{i}}_{\ff}   \)  \sum_{\cup_{i=1}^{m} C_{i}=A^{c}}\,\,\,\prod_{i=1}^{m}\al \,\mbox{cy}(C_{i}).\nonumber
 \end{eqnarray}
This follows from (\ref{pe.16}), which says the $Q^{x_{0},x_{0}}\(\prod_{i\in A}   L^{ x_{i}}_{\ff}   \)=\mbox{cy}(A\cup \{0\})$.\qed

Comparing (\ref{pp.15r}) with (\ref{gp.17}) proves that in symmetric case,  with $\al=1/2$, we have 
\begin{equation}
\{ \wh L_{\al}^{x}, x\in S\}\overset{law}{=} \{  {1 \over 2}G^{2}_{ x}, x\in S\}   \label{pp.11}
\end{equation}
In this case,  (\ref{pp.9v}) is Dynkin's isomorphism theorem
(\ref{dit.1}) which we write as 
  \begin{equation}
   E_{G} \(G^{2}_{ x}\,\,  F\(   {1 \over 2}G^{2}_{ x_{i}} \) \)=E_{G}Q^{x,x}\(F\(   {1 \over 2}G^{2}_{ x_{i}}+L^{ x_{i}}_{\ff} \) \).
   \label{pp.10}
\end{equation}
(\ref{pp.11}) explains why, intuitively, 
Gaussian squares should  be related to Markov local times. The reason is that Gaussian squares are themselves sums of local times.

 \subsection{A Palm formula proof of the isomorphism theorem}

 In this section we show that  the isomorphism theorem (\ref{pp.9v}) is just a simple application of the Palm formula.
We   apply the Palm formula (\ref{pp.4}) with $f(\om)=  L_{\ff}^{x}(\om)$ and $G(\mathcal{L})=F\( \wh L_{\al}^{x_{j}}\)$ where as before
 \begin{equation}
 \wh L_{\al}^{x}=N(L_{\ff}^{x})=\sum_{\om\in \mathcal{L_{\al}}} L_{\ff}^{x}(\om).
 \label{pp.5}
\end{equation}
Using the fact that $\mu$ is non-atomic, we see that for  any fixed $\om'\in \Om_{\De}$, almost surely  $\om' \notin \mathcal{L}_{\al}$, and consequently
 \begin{equation}
\wh L_{\al}^{x_{j}}(\om'\cup\mathcal{L_{\al}})
=\sum_{\om\in \mathcal{\om' \cup L_{\al}}}  L_{\ff}^{x_{j}} (\om)=\wh L_{\al}^{x_{j}}(\mathcal{L_{\al}})+L_{\ff}^{x_{j}}(\om').\label{pp.6}
\end{equation}
 Thus  
\begin{equation}
G(\om'\cup\mathcal{L}_{\al})=F\(\wh L_{\al}^{x_{j}}(\mathcal{L_{\al}})+L_{\ff}^{x_{j}}(\om')\).\label{pp.7}
\end{equation}  
Then by the Palm formula (\ref{pp.4}) 
 \begin{equation}
E_{\mathcal{L}_{\al}} \(\wh L_{\al}^{x} \,F\(\wh L^{x_{j}}\)  \)
=\al E_{\mathcal{L}_{\al}}  \int \(L_{\ff}^{x}(\om') \,  F \( \wh L_{\al}^{x_{j}} +L_{\ff}^{x_{j}}(\om')\)\,d\mu (\om')\). \label{pp.8}
\end{equation}
 It follows from (\ref{pe.17}) that we can rewrite this as 
 \begin{equation}
E_{\mathcal{L}_{\al}} \(\wh L_{\al}^{x} \,F\(\wh L_{\al}^{x_{j}}\)  \)
=\al E_{\mathcal{L}_{\al}} Q^{x,x}\(  F \( \wh L_{\al}^{x_{j}} +L_{\ff}^{x_{j}}\) \), \label{pp.9}
\end{equation}
which is  (\ref{pp.9v}).

\subsection{Permanental processes}

Our goal in this sub-section is to better understand the stochastic process $\{ \wh L_{\al}^{x}, x\in S\}$ which appears in our isomorphism theorem (\ref{pp.9}).

Using  (\ref{pp.15r}), (\ref{pe.loop})   and writing  $B_{i}=\{ x_{i,1}, x_{i,2},\ldots, x_{i,|B_{i}|}\}$ we have 
\be
E_{\mathcal{L}_{\alpha}}\( \prod_{j=1}^{n}\wh L_{\al}^{x_{j}}  \)=\sum_{\ell=1} ^{n}\,\,\sum_{\cup_{i=1}^{\ell} B_{i}=[1,n]}\al^{\ell}\,\prod_{i=1}^{\ell}\,\sum_{\pi\in \mathcal{P}^{\odot}_{|B_{i}|}}\,u(x_{i,\pi_{1}},x_{i,\pi_{2}})\cdots u(x_{i,\pi_{|B_{i}|}},x_{i,\pi_{1}}).\label{pp.16}
\ee
(\ref{pp.16}) can also be written as 
\be
E_{\mathcal{L}_{\alpha}}\( \prod_{j=1}^{n}\wh L_{\al}^{x_{j}}  \)
=\sum_{\pi\in\mathcal{P}_{n}} \al^{c(\pi)} \,u(x_{1},x_{\pi_{1}})u(x_{2},x_{\pi_{2}})\cdots u(x_{n},x_{\pi_{n}}).\label{pp.17} 
\ee
In particular
\begin{equation}
E_{\mathcal{L}_{\alpha}}\(\wh L_{\al}^{x}  \)=\al u( x,x),\hspace{.2in}\mbox{Cov}\(\wh L_{\al}^{x}, \wh L_{\al}^{y}  \)
=\al u( x,y)u( y,x).\label{covz}
\end{equation}
 When $\al=1$ the right hand side of  (\ref{pp.17}) is the permanent of the matrix $\{   u( x_{i}, x_{j})\}$, while if $\al=-1$ we obtain the determinant. 
 In general, this is referred to as the $\al$-permanent, see \cite{VJ}, and a process satisfying (\ref{pp.17}) is called an $\al$-permanental process.
 
 By (\ref{covz}), $u( x,y)u( y,x)$ is positive definite, hence so is  $\sqrt{u( x,y)u( y,x)}$. Let $( G_{x}, G_{y})$ be the Gaussian random vector with covariance $\sqrt{u( x,y)u( y,x)}$.
 An important property of the  $1/2$-permanental process $\wh L_{\al}^{x},\, x\in S$ is that the  bivariate distributions $(\wh L_{\al}^{x}, \wh L_{\al}^{y} )$ are the same as  $( G_{x}^{ 2}/2, G_{y}^{ 2}/2)$. To see this, it suffices to show that     
 \begin{equation}
 E_{\mathcal{L}_{\alpha}}\(\left(  \wh L_{\al}^{x} \right)^{ j} \left(  \wh L_{\al}^{y} \right)^{k}\)=E\left(  \left( G_{x}^{ 2}/2 \right)^{ j} \left(G_{y}^{ 2}/2\right)^{k} \right)\label{covx}
 \end{equation}
   for all $j,k$. 
 Comparing  (\ref{pp.15r}) with (\ref{gp.17})   with $\al=1/2$, shows that both involve cycles, the only difference being that the left hand side involves cycles with respect to $\sqrt{u( x,y)u( y,x)}$  while the right hand side uses $u( x,y)$.  $u( x,y)$ is the same as $\sqrt{u( x,y)u( y,x)}$ when $x=y$, but note that in the left hand side of (\ref{covx}), whenever we have $u( x,y)$ with $x\neq y$, we must also have a corresponding $u( y,x)$. If we replace both elements of this pair by $\sqrt{u( x,y)u( y,x)}$ we will not change the value of the left hand side. Implementing this change for all $u( x,y)$ with $x\neq y$ establishes (\ref{covx}).
 
 The importance of the fact that $(\wh L_{\al}^{x}, \wh L_{\al}^{y} ) \overset{dist}{=} ( G_{x}^{ 2}/2, G_{y}^{ 2}/2)$ comes from the fact that in proving the sufficiency of the condition (\ref{dit.12}) for the continuity of Gaussian processes, all that is used is the bivariate distributions. This allows us to obtain a similar result for permanental processes, see \cite{MRperm}.
 
 {\it   See \cite{Le Jan1} and the earlier Arxiv version of \cite{MRperm}. See \cite{FPY} for Markovian bridges.   For Poisson processes see \cite{K}. The Palm formula is given in \cite[Lemma 2.3]{Bertoin}. This reference assumes that $S$ is Polish, but that assumption is not necessary. Permanental processes were introduced in \cite{VJ}, and their relevance to isomorphism theorems was established in \cite{EK}. For later developments see \cite{FR, LMR, LMR2}.  For other work on loop soups see \cite{LL, LF, LW}.  }

\section{A Poisson process approach to the generalized second Ray-Knight theorem}

Using excursion theory, we can give a simple proof of the generalized  second Ray--Knight Theorem   which does not make require us to work with $\tau(\la)$ for an independent exponential $\la$.

As before we assume that $X$ is symmetric, recurrent, with $P^{x}\(T_{0}<\ff\)=P^{x}\(T_{0}<\ff\)=1$ for all $x\in S $ and $u(0,0)=\ff$. We let ${\bf  n}$ denote the excursion measure for $X$ with respect to the point $0$. ${\bf  n}$ is a $\si$-finite measure on   $\Om_{\De} $. Let $\mathcal{E}_{t}$ be  a Poisson   process  on $\Om_{\De}$ with intensity measure $t {\bf  n}$. It is a fundamental result of 
excursion theory, \cite{B, DM3} that  
 \begin{equation}
\{   L_{\tau(t)}^{x}, x\neq 0, P^{ 0}\}=\{   N(L_{\ff}^{x})=\sum_{\om\in \mathcal{E}_{t} } L_{\ff}^{x}(\om), x\neq 0, P_{\mathcal{E}_{t}}\}.
 \label{et.1}
\ee
Hence by the moment formula (\ref{pp.15})
\begin{equation}
P^{ 0}\( \prod_{j=1}^{n}L_{\tau(t)}^{x_{j}} \)=E_{ \mathcal{E}_{t} }\( \prod_{j=1}^{n}N(L_{\ff}^{x_{j}}) \)=\sum_{\cup_{i=1}^{\ell} B_{i}=[1,n]} \,\,\prod_{i=1}^{\ell}\,t \,{\bf  n}\(\prod_{j\in B_{i}} L_{\ff}^{x_{j}}\).\label{}
\end{equation}
In view of (\ref{rk.5}) we need only show that
\begin{equation}
{\bf  n}\(\prod_{j\in B_{i}} L_{\ff}^{x_{j}}\)=\mbox{ ch}_{0} (B_{j}). \label{et.111}
\end{equation} 
This is the content of the next Lemma.

\subsection{Excursion local time} \label{sec-exlt}

\bl \label{lem-ltex} For any $y\neq 0$
\begin{equation}
{\bf  n}\(L^{y}_{\ff}\)=1,\label{39.1}
\end{equation}
and for $y_{1},\cdots,y_{k}\neq 0$, $k\geq 2$, 
\be
{\bf  n}\(\prod_{j=1}^{k}L_{\ff}^{y_{j}}\) = \sum_{\pi\in \mathcal{P}_{k}} \prod_{j=1}^{k-1} u_{T_0}(y_{\pi(j)},y_{\pi(j+1)}).\label{39.2}
\ee
\el

{\bf  Proof:  }Under ${\bf  n}$, the coordinate process is Markovian with an entrance law which we denote by $\iota_{t},\,t>0,$ and transition probabilities given by the stopped process $X^{0}_{t}=X_{t\wedge T_{0}}$.  $X^{0}$ has potential densities $u_{T_{0}}(x,y)$ for $x,y\neq 0$. where $u_{T_{0}}(x,y)$ are the potential densities for the process obtained by killing $X$ at $T_{0}$.

Let $\iota_{t},\,t>0,$ denote the entrance law for ${\bf  n}$. Then
\begin{equation}
{\bf  n}\(f\(X_{t})\)\)=\iota_{t}(f).\label{39.3}
\end{equation}
It follows from \cite[XV,  (78.3)]{DM3}  or \cite[VI, (50.3)]{RW} that for any  $f$ which is zero at $0$
\begin{equation}
\int_{0}^{\ff}e^{-\al t}\iota_{t}(f)\,dt={1 \over u^{\al}(0,0)}\int u^{\al}(0,x) f(x)\,dm(x).\label{39.4}
\end{equation}
Hence, if $f_{y,\ep}$ is an approximate $\de$-function for $y$ supported in the ball of radius $\ep$ centered at $y$, then for  $\ep $ sufficiently small
\begin{equation}
{\bf  n}\(\int_{0}^{\ff} e^{-\al t}f_{y,\ep}\(X_{t}\)\,dt\)={1 \over u^{\al}(0,0)}\int u^{\al}(0,x) f_{y,\ep}(x)\,dm(x).\label{39.5}
\end{equation}
We claim that for any bounded measurable function $g$
\begin{equation}
\lim_{\ep\to 0}\int g( t)f_{y,\ep}\(X_{t}\)\,dt=\int g( t)\,dL_{t}^{y}.\label{39.5ap}
\end{equation}
It suffices to prove this for $g$ of the form $g( t)=1_{\{   [0,r]\}}( t)$, in which case it follows from (\ref{mp.6}).
Hence, letting $\ep\to 0$ in (\ref{39.5}) and then using (\ref{rk.21}) gives 
\begin{equation}
{\bf  n}\(\int_{0}^{\ff} e^{-\al t}\,dL_{t}^{y}\)={u^{\al}(0,y) \over u^{\al}(0,0)}= E^{y}\(e^{-\al T_{0}}\).\label{39.6}
\end{equation}
Letting $\al\to 0$ gives 
\begin{equation}
{\bf  n}\(L^{y}_{\ff}\)=P^{y}\(T_{0}<\ff\)=1\label{39.6a}
\end{equation}
by our assumption. This proves (\ref{39.1}).

Now let  $0< t_{1}\leq \cdots \leq t_{k-1}\leq t_{k} < \ff$. Then, if $P^{0}_{t}$ denotes the transition operator for the stopped process $X^{0}_{t}$,
\begin{eqnarray}
&&{\bf  n}\(\prod_{j=1}^{k}  f_{y_{j},\ep}\(X_{t_{j}}\) \)
\label{39.7}\\
&&=\int f_{y_{1},\ep}(x_{1})
\prod_{j=2}^{k}P^{0}_{t_{j}-t_{j-1}}(x_{j-1},dx_{j})f_{y_{j},\ep}(x_{j})\,\iota_{t_{1}}(dx_{1}),\nn   \nonumber
\end{eqnarray}
Hence, using  (\ref{39.4}) with $\ep$ sufficiently small
\begin{eqnarray}
&&{\bf  n}\(\int_{\{0< t_{1}\leq \cdots \leq t_{k-1}\leq t_{k} < \ff\}}e^{-\al t_{k}}\prod_{j=1}^{k}  f_{y_{j},\ep}\(X_{t_{j}}\) \,dt_{j}\)
\label{39.8}\\
&&={1 \over u^{\al}(0,0)}\int u^{\al}(0,x_{1})  f_{y_{1},\ep}(x_{1})
\prod_{j=2}^{k}u^{\al}_{T_0}(x_{j-1},x_{j})f_{y_{j},\ep}(x_{j})\prod_{j=1}^{k} \, dm(x_{j}).\nn   \nonumber
\end{eqnarray}
Letting $\ep\to 0$ gives 
\begin{eqnarray}
&&{\bf  n}\(\int_{\{0< t_{1}\leq \cdots \leq t_{k-1}\leq t_{k} < \ff\}}e^{-\al t_{k}}\prod_{j=1}^{k}\,dL_{t_{j}}^{y_{j}} \)
\label{39.8}\\
&&\hspace{2 in}={u^{\al}(0,y_{1})   \over u^{\al}(0,0)}  
\prod_{j=2}^{k}u^{\al}_{T_0}(y_{j-1},y_{j}),\nn   \nonumber
\end{eqnarray}
and (\ref{39.2}) follows as before on letting $\al\to 0$.\qed

{\it The passage to the limit  $\ep\to 0$ under the measure ${\bf  n}$ in both (\ref{39.6}) and (\ref{39.8}) requires more justification. However,  ${\bf  n}(1_{\{\ze\geq \de\}}\cdot )$ is a finite measure for any $\de>0$. Using the material we have presented, it is easy to check that the integrands are uniformly bounded in $L^{2}(\,d{\bf  n})$  and hence uniformly integrable in $L^{2}(1_{\{\ze\geq \de\}}\,d{\bf  n})$. We can thus take the $\ep\to 0$ in $L^{2}(1_{\{\ze\geq \de\}}\,d{\bf  n})$, and then the  $\de\to 0$ limit using the monotone convergence theorem.}

\section{Another  Poisson process isomorphism theorem: random interlacements}

Sznitman has recently developed an isomorphism theorem related to a Poisson process for quasi-processes which he refers to as random interlacements, \cite{Sz2}. This  isomorphism theorem has the structure of the generalized second Ray-Knight theorem. Typically, the underlying Markov processes (think of Brownian motion in two or more dimensions) do not have finite potential densities, hence there are no local times nor associated Gaussian processes $G_{x}$ indexed by points in the state space $S$. We first develop material for the associated Gaussian process which is now  indexed by measures on $S$. We then introduce quasi-processes and random interlacements. At that stage the isomorphism theorem will be straightforward.

\subsection{Gaussian fields}

We assume that $X$ is a symmetric Markov process in $S$ with transition densities $p_{t}(x,y)$. As before, these are positive definite and consequently the potential densities $u(x,y)$ will be 
positive definite in the wide sense, that is
\begin{equation}
\int\int u(x,y)\,d\nu(x)\,d\nu(y)\geq 0\label{gf.1}
\end{equation}
for any positive measure $\nu$ on $S$.
Let $\mathcal{M}$ denote the set of finite positive measures on $S$ and
\begin{equation}
 \mathcal{G}^{1}=\{\nu\in \mathcal{M}\,|\, \int\int u(x,y)\,d\nu(x)\,d\nu(y)<\ff\}.\label{q2.2axc}
\end{equation}
Let $G_{\nu}$ denote the mean zero Gaussian process on $\mathcal{G}^{1}$ with covariance
\begin{equation}
E\( G_{\nu}G_{\nu'}\)=\int\int u(x,y)\,d\nu(x)\,d\nu'(y).\label{q2.2}
\end{equation}

We would like to find an analogue of $\int G_{x}^{2}\,d\nu(x)$ to obtain some version of (\ref{gp.17}), but if any of the sets $A_{l}$ in (\ref{gp.17}) are singletons, then the cycle term would be $\int u(x,x)\,d\nu(x)$ and for the processes we would like to consider, $u(x,x)=\ff$ for all $x$. It is the need to eliminate such singletons that leads us to define the Wick square.
In the following we assume that $u_{\de}(x,y)=:\int_{\de}^{\ff}p_{t}(x,y)\,dt<\ff$ for all $x,y\in S$, and $\de>0$.

It is then easy to check that $p_{\ep}(x,y)\,dm(y)\in \mathcal{G}^{1}$ for any $\ep>0$ and $x\in S$. Set
\begin{equation}
G_{x,\ep}=G_{p_{\ep}(x,y)\,dm(y)}.\label{gf.2}
\end{equation}
Then $G_{x,\ep}$ is a Gaussian process on $S\times (0,\ff)$ with covariance
\begin{equation}
E\( G_{x,\ep}G_{x',\ep'}\)=  u_{\ep+\ep'}(x,y).\label{gf.3}
\end{equation}
If we set $:G_{x,\ep}^{2}:=G_{x,\ep}^{2}-E\(G_{x,\ep}^{2}\)$ then it is easy to see that 
\begin{equation}
E\(\prod_{i=1}^{n}:G^{2}_{x_{i},\ep}:/2\)=\sum_{\stackrel{A_{1}\cup \cdots\cup A_{j}=[1,n]}{|A_{l}|\geq 2}}\prod^{j}_{l=1}{1 \over 2}\,\mbox{cy}_{\ep}(A_{l}),\label{gf.4}
\end{equation}
where
\begin{equation}
\mbox{cy}_{\ep}(A_{l})=\sum_{\pi\in \mathcal{P}^{\odot}_{| A_{l}|}}u_{2\ep}(x_{l_{\pi(1)}},  x_{l_{\pi(2)}})\cdots u_{2\ep}(x_{l_{\pi (| A_{l}|)}},  x_{l_{\pi(1)}}).\label{gf.5}
\end{equation}

Define the  the Wick square   
    \begin{equation}
  :G^{2}: (\nu)=\lim_{\ep\rar0} \int :G^{2}_{x,\ep}: \,d\nu(x). \label{a.3a}
    \end{equation}
Using (\ref{gf.3}) we can showk that if 
\begin{equation}
\int\int u^{2}(x,y)\,d\nu(x)\,d\nu(y)<\ff,\label{a.3a4}
\end{equation}
then   the limit in (\ref{a.3a}) exists in all $L^{p}$, and we have
\begin{equation}
E\(\prod_{i=1}^{n}:G^{2}: (\nu_{i})/2\)=\sum_{\stackrel{A_{1}\cup \cdots\cup A_{j}=[1,n]}{|A_{l}|\geq 2}}\prod^{j}_{l=1}{1 \over 2}\,\mbox{cy}(A_{l},\nu),\label{gf.4}
\end{equation}
where
\begin{equation}
\mbox{cy}(A_{l},\nu)=\sum_{\pi\in \mathcal{P}^{\odot}_{| A_{l}|}}\int u(x_{l_{\pi(1)}},  x_{l_{\pi(2)}})\cdots u(x_{l_{\pi (| A_{l}|)}},  x_{l_{\pi(1)}})\prod_{j\in A_{l}}\,d\nu_{j}(x_{j}).\label{gf.5}
\end{equation}
 (See \cite[Lemma 3.3]{MRcontperm} for an important ingredient in the proof ).  
 Set
\begin{equation}
 \mathcal{G}^{2}=\{\nu\,|\, \int\int u^{2}(x,y)\,d\nu(x)\,d\nu(y)<\ff\},\label{q2.2a}
\end{equation}
and let $\mathcal{G}_{K}^{2}$ denote the subset of  measures $\nu\in \mathcal{G}^{2}$ with  support in the compact set $K\subseteq S$.
 
Let $|\nu|$ denote the mass of $\mu$. Exactly as in (\ref{dit.9}) we can then show that for $\nu_{i}\in \mathcal{G}^{2}$
\bea
&&
E_{G} \( \prod_{i=1}^{k}\( :G^{2}: (\nu_{i})/2+\sqrt{2t}\,G_{\nu_{i}} +t|\nu_{i}|\)   \)\label{gf.6}\\
&&=\sum_{\stackrel{\cup_{i=1}^{l} A_{i}\cup_{j=1}^{m} B_{j}=[1,k]}{|A_{l}|\geq 2}}\,\,
\prod_{i=1}^{l}{1 \over 2} \mbox{cy} (A_{i},\nu)    \prod_{j=1}^{m}t\,\,\mbox{ch} (B_{j},\nu)
\nn
\eea
where $\mbox{ch}(B)=|\nu_{i}|$  if $B=\{i\}$ and, if  $|B|>1$ with $B=\{b_{1},b_{2},\cdots, b_{|B|}\}$  then the chain function $\mbox{ch}(B,\nu)$ is defined as 
\begin{equation}
\mbox{ch}(B,\nu)=\sum_{\pi\in \mathcal{P}_{| B|}}\int u(x_{b_{\pi(1)}},  x_{b_{\pi(2)}})\cdots u(x_{b_{\pi (| B|-1)}},  x_{b_{\pi (| B|)}})\prod_{j\in B}\,d\nu_{j}(x_{j}).\label{gf.7}
\end{equation} 

It follows as in proof of the generalized second  Ray-Knight theorem that if we can find a family of random variables $\{  S_{\nu,t}, \nu\in  \mathcal{G}^{2} \}$ such that  
\be
 P\(   \prod_{i=1}^{k}S_{\nu_{i},t} \) 
  = \sum_{m=1}^{k}\sum_{\stackrel{\mbox{\scriptsize unordered}}{B_{1}\cup\cdots\cup B_{m}=[1,k]}} \,\,t^{m}\prod_{j=1}^{m} \mbox{ ch} (B_{j},\nu), \label{gf.8}
 \ee
then we will have established the isomorphism theorem
\begin{equation}
E_{G}P\(F\(  S_{\nu_{i}}+:G^{2}: (\nu_{i})/2\) \)= E_{G}\(F\(:G^{2}: (\nu_{i})/2+\sqrt{2t}\,G_{\nu_{i}} +t|\nu_{i}|\ \) \),\label{gf.9}
\end{equation}
Such random variables $S_{\nu}$ will come from additive functionals of  random interlacements.

\subsection{Quasi-processes and additive functionals}\label{sec-qp}

 Let 
$X=\left(\Omega,\mathcal{F},\mathcal{F}_t,X_t, \theta_t,P^x
\right)$ be a `nice'  symmetric  transient Markov process as before  with LCCB state space $S$ and  transition
densities $p_{t}(x,y)$ with respect to a $\si$-finite measure $m$. We assume that $m$ is dissipative, that is, that $\int u(x,y)f(y)\,m(dy)<\ff$ $m$-a,e for each non-negative  $f\in L^{1}(m)$. This will hold for example for Brownian motion in $R^{ 3}$ or exponentially killed Brownian motion in $R^{ 2}$, with $m$ being Lebesgue measure.

Let $W$ denote the set of paths $\om: R^{1}\mapsto S\cup \De$ which are $S$ valued and right continuous on some open interval $(\al (\om), \bb (\om))$ and $\om (t)=\De $ otherwise. Let $Y_{t}=\om (t)$,   and define the shift operators
\begin{equation}
(\si_{t}\om) (s)=\om (t+s),\hspace{.2 in}s,t\in R^{1}.\label{q.1}
\end{equation}
Set $\mathcal{H}=\si\(Y_{s}, s\in R^{1}\)$ and $\mathcal{H}_{t}=\si\(Y_{s}, s\leq t\)$. 
Let $\mathcal{A}$ denote the $\si$-algebra of shift invariant events in $\mathcal{H}$. 
 The quasi-process associated with $X$ is the measure $\bf{P}_{m}$ on $\(W,\mathcal{A} \)$ which satisfies the following two conditions: 
 \begin{equation}
\hspace{-.5in}(i):\hspace{.5in}\mbox{$ \bf{P}_{m}$}\(\int_{R^{1}}f\(Y_{t}\)\,dt\)=m(f),\label{q.2}
 \end{equation}
 for all measurable $f$ on $S$,
 and (ii):  for any intrinsic stopping time $T$, $Y_{T+t}, t>0$ is Markovian with semigroup $P_{t}$, recall (\ref{mp.1a}), under 
 $\mbox{$ \bf{P}_{m}$}|_{\{T\in R\}}$.
 An $\mathcal{H}_{t^{+}}$ stopping time $T$ is called intrinsic if $\al\leq T\leq\bb$ on $\{T<\ff\}$ and 
 $T=t+T\circ \si_{t}$ for all $t\in R^{1}$. A first hitting time is an example of an intrinsic stopping time.
 
 If $L^{\nu}_{t}, t\geq 0$ denotes the continuous additive functional, (recall (\ref{mp.7})), on $\Omega$ with  $E^{ x}\left(  L^{\nu}_{\ff} \right)= \sup_{x}\int u(x,y)\,d\nu(y)<\ff$ then there is an extension to $W$, which we also denote by $L^{\nu}_{t}, t\in R^{1}$  with the property that
  \begin{equation}
\mbox{$ \bf{P}_{m}$}\(\int_{R^{1}}g\(Y_{t}\)\,dL^{\nu}_{t}\)=\nu (g),\label{q.3}
 \end{equation}
for all measurable $g$, see \cite[XIX, (26.5)]{DM4}. For example, if $\nu=f\,dm$ then
\begin{equation}
L^{f\,dm}_{t}=\int_{-\ff}^{t}f\(Y_{s}\)\,ds\label{q.3a}
\end{equation}
and (\ref{q.3}) follows easily from (\ref{q.2}).
In general one can think of $L^{\nu}_{t}$ as
\begin{equation}
L^{\nu}_{t}=\lim_{\ep\to 0}\int_{S} \int_{-\ff}^{t}f_{x,\ep}\(Y_{s}\)\,ds\,d\nu(x).\label{q.3b}
\end{equation}

 \bl \label{lem-qpm} For any $\nu_{1},\cdots,\nu_{k}$, with support in some compact $K\subset S$
\be
\mbox{$ \bf{P}_{m}$}\(\prod_{j=1}^{k}L_{\ff}^{\nu_{j}}\) = \sum_{\pi\in \mathcal{P}_{k}} \int \prod_{j=1}^{k-1} u(y_{j},y_{j+1}) \prod_{j=1}^{k}\,d\nu_{\pi (j)}(y)=\mbox{ch } ([1,k],\nu). \label{q.4}
\ee
\el

Compare Lemma \ref{lem-ltex}.

{\bf  Proof:  } Let  $T_{K}$ denote the first hitting time of $K$. Since the measures $\nu_{i}$ are supported in $K$, it follows that the functionals $L_{t}^{\nu_{i}}$ do not grow until time $T_{K}$. Hence
\begin{eqnarray}
&&\mbox{$ \bf{P}_{m}$}\(\int_{\{-\ff<t_{1}\leq \cdots \leq t_{k-1}\leq t_{k} < \ff\}}\prod_{j=1}^{k}  \,dL_{t_{j}}^{\nu_{j}}\)
\label{q.5}\\
&&=\mbox{$ \bf{P}_{m}$}\(\int_{\{0\leq t_{1}\leq \cdots \leq t_{k-1}\leq t_{k} < \ff\}}\prod_{j=1}^{k}  \,dL_{T_{K}+t_{j}}^{\nu_{j}}\).\nn   \nonumber
\end{eqnarray}
Hence by the second property of $\mbox{$ \bf{P}_{m}$}$
this equals
\begin{equation}
\mbox{$ \bf{P}_{m}$}\(\int_{0}^{\ff} h\(Y_{T_{K}+t_{1}}\)\,dL_{T_{K}+t_{1}}^{\nu_{1}}\),\label{q.6}
\end{equation} 
where
\bea
h(x)&=&E^{x}\(\int_{\{0\leq t_{2}\leq \cdots \leq t_{k-1}\leq t_{k} < \ff\}}\prod_{j=2}^{k}  \,dL_{t_{j}}^{\nu_{j}}\)\label{q.7}\\
&=&\int u(x,y_{2}) \prod_{j=2}^{k-1} u(y_{j},y_{j+1}) \prod_{j=2}^{k}\,d\nu_{j}(y).\nn
\eea
(For those unfamiliar with such calculations, think of (\ref{q.3a}) or more generally (\ref{q.3b})).
Using  once again  the fact that $L_{t}^{\nu_{1}}$ doesn't  grow until time $T_{K}$ and then (\ref{q.3})   shows that
\begin{eqnarray}
&&\mbox{$ \bf{P}_{m}$}\(\int_{\{-\ff<t_{1}\leq \cdots \leq t_{k-1}\leq t_{k} < \ff\}}\prod_{j=1}^{k}  \,dL_{t_{j}}^{\nu_{j}}\)
\label{q.5}\\
&&=\mbox{$ \bf{P}_{m}$}\(\int_{R^{1}}h\(Y_{t_{1}}\)\,dL^{\nu_{1}}_{t_{1}}\)\nn\\
&&= \int \prod_{j=1}^{k-1} u(y_{j},y_{j+1}) \prod_{j=1}^{k}\,d\nu_{j}(y),\nn   \nonumber
\end{eqnarray}
and (\ref{q.4}) follows.
\qed

\subsection{Interlacements}\label{sec-intsym}

Interlacements are the soup of a quasi-process. More precisely, the interlacement $\mathcal{I}_{t}$ is the Poisson process with intensity measure \mbox{$t  \bf{P}_{m}$}. We let 
${P}_{\mathcal{I}_{t}}$ denote probabilities for the process $\mathcal{I}_{t}$. Let
\begin{equation}
\wt L^\nu_\ff=\sum_{\om \in \mathcal{I}_{t}}L^\nu_\ff(\om).\label{q2.1}
\end{equation}

Using (\ref{q.4}) and the moment formula (\ref{}), we see that the functionals  $ \wt L^\nu_\ff$, under the measure 
${P}_{\mathcal{I}_{t}}$, satisfy (\ref{gf.8}). In view of (\ref{gf.9}) we have the following interlacement Isomorphism theorem which is essentially due to Sznitman, \cite{Sz2}.

\bt
For any  $t >0$, compact $K\subset S$ and countable  $D\subseteq  \mathcal{G}_{K}^{2}$,
\bea
&&
\Big\{ \wt L^\nu_\ff+\textstyle{ 1\over 2}:G^{2}:(\nu),\,\nu\in D, {P}_{\mathcal{I}_{t}}\times P_{G}\Big\}\label{q2.0}\\
&&\stackrel{law}{=}
\Big\{\textstyle{ 1\over 2}:G^{2}:(\nu)+\sqrt{2t}G_{\nu}+t  |\nu|,\,\nu\in D,P_{G} \Big\}.\nn
\eea
\et

\section{Isomorphism theorems via Laplace transforms}

In this section we give alternate proofs for our   Isomorphism theorems. The innovation here is that we use the moment generating function   of Gaussian squares, described in the next subsection, instead of the Gaussian moment formulas of Section \ref{sec-gmf}. On the other hand, we still need the local time moment formulas, and in particular for the generalized second Ray-Knight theorem we have seen that the derivation is not trivial.  

\subsection{Moment generating functions  of Gaussian squares}\label{sec-gmgf}

Let   $G=(G_{1}, \ldots, G_{n})\in R^{n}$ be a Gaussian random vector with covariance matrix $C$. 
If  $C$ is invertible we first show that for all bounded measurable functions $F$ on $R^{d}$
\begin{equation}
E\(F(G_{1}, \ldots, G_{n})\)={1 \over (2\pi)^{n/2}\sqrt{|C|}}\int_{R^{n}} F(x)e^{-(x,C^{-1}x)/2}\,dx\label{lt.1}
\end{equation}
where $|C|$ denotes the determinant of $C$. To see this it suffices to prove it for $F$ of the form 
$F(x)=e^{i(y,x)}$, in which case we need to show that
\begin{equation}
E\(e^{i(y,G)}\)={1 \over (2\pi)^{n/2}\sqrt{|C|}}\int_{R^{n}} e^{i(y,x)}e^{-(x,C^{-1}x)/2}\,dx.\label{lt.2}
\end{equation}
Setting $x=C^{1/2}z$, (recall the paragraph following (\ref{gp.8})), so that $dx=|C|^{1/2}\,dz$, the right hand side of (\ref{lt.2}) becomes 
\begin{equation}
{1 \over (2\pi)^{n/2}}\int_{R^{n}} e^{i(C^{1/2}y,z)}e^{-(z, z)/2}\,dz=e^{-(y,Cy)/2}\label{lt.3}
\end{equation}
which, by (\ref{gp.7}) equals the left hand side of  (\ref{lt.2}).

 We now show that for any Gaussian random vector $G=(G_{1}, \ldots, G_{n})$ with covariance matrix $C$, any vector $u=(u_{1}, \ldots, u_{n})$ and $\la_{1}, \ldots, \la_{n}$ sufficiently small
 \begin{equation}
 E\(e^{\sum_{j=1}^{n}\la_{j}u_{j}G_{j}+ \la_{j}G^{2}_{j}/2    }\)={1 \over  \sqrt{|I-\La C|}}
e^{(u,\La \bar C\La u)/2}\label{lt.4}
 \end{equation}
where $\La$ is the diagonal matrix with entries $(\la_{1}, \ldots, \la_{n})$ and
\begin{equation}
\bar C=(I- C\La)^{-1}C.\label{lt.5}
\end{equation}
(\ref{lt.4}) will be the key to the alternative proofs of the Isomorphism theorems given in this section.

Proof of (\ref{lt.4}): Assume first that $C$ is invertible. Then from (\ref{lt.5})
\begin{equation}
\bar C^{-1}=C^{-1}(I- C\La)=C^{-1}-\La.\label{lt.5a}
\end{equation}
Hence, using (\ref{lt.1}) we have
\begin{eqnarray}
&& E\(e^{\sum_{j=1}^{n}\la_{j}u_{j}G_{j}+ \la_{j}G^{2}_{j}/2    }\)\label{lt.6}\\
&&={1 \over (2\pi)^{n/2}\sqrt{|C|}}\int_{R^{n}} e^{(\La u, x)}e^{(x,\La x)/2}e^{-(x,C^{-1}x)/2}\,dx
   \nonumber\\
&&={1 \over (2\pi)^{n/2}\sqrt{|C|}}\int_{R^{n}} e^{(\La u,x)} e^{-(x,\bar C^{-1}x)/2}\,dx
   \nonumber
\end{eqnarray}
It is clear from (\ref{lt.5a}) that for  $\la_{1}, \ldots, \la_{n}$ sufficiently small, $\bar C^{-1}$ is invertible, symmetric and positive definite. Hence changing variables $x=\bar C^{1/2}z$ as before the last display
\begin{equation}
={\sqrt{|\bar C|} \over \sqrt{|C|}}{1 \over (2\pi)^{n/2}}\int_{R^{n}} e^{(\bar C^{1/2}\La  u,z)} e^{-(z,z)/2}\,dz=\sqrt{  |\bar C|  \over  |C|}e^{(u,\La \bar C\La u)/2},\label{lt.7}
\end{equation}
and (\ref{lt.4}) for $C$   invertible follows from (\ref{lt.5}).

For general $C$, recall that we can find an orthonormal system of vectors $u_{i}, 1\leq i\leq n$ such that $Cu_{i}=c_{i}u_{i}, 1\leq i\leq n$, and the fact that $C$ is   positive definite implies that all $c_{i}\geq 0$. We can then define the matrix $C_{\ep}$ by setting $C_{\ep}u_{i}=(c_{i}+\ep)u_{i}, 1\leq i\leq n$. $C_{\ep}$ is clearly symmetric, positive definite  and invertible. We then obtain (\ref{lt.4}) by first proving it for the Gaussian random vector $G_{\ep}$  with covariance matrix $C_{\ep}$ and then taking the limit as $\ep\rar 0$. That $G_{\ep}\rar G$ in distribution follows from (\ref{gp.7}).
\qed

We note that (\ref{lt.4}) immediately implies that
 \begin{equation}
 E\(e^{\sum_{j=1}^{n}\la_{j}(G_{j}+ u_{j})^{2}/2    }\)={1 \over  \sqrt{|I-\La C|}}e^{(u,\La u)/2}
e^{(u,\La \bar C\La u)/2}.\label{lt.8}
 \end{equation}
We also note the following computation for later use:
\begin{eqnarray}
\La+\La \bar C\La&=&\La+\La(I- C\La)^{-1}C\La
\label{lt.9}\\
&=&\La\(I+(I- C\La)^{-1}C\La\)=\La\(\sum^{\ff}_{k=0} (C\La)^{k} \).   \nonumber
\end{eqnarray}

\subsection{Another proof of the Dynkin Isomorphism theorem}

It suffices to show that 
\begin{equation}
 E_{G}Q^{x_{1},x_{2}}\(e^{\sum_{j=1}^{n}\la_{j}(L^{ x_{j}}_{\ff}+{1 \over 2}G^{2}_{ x_{j}})}  \)=  E_{G} \(G_{ x_{1}}G_{ x_{2}}\, \exp^{\sum_{j=1}^{n}\la_{j}{1 \over 2}G^{2}_{ x_{j}}} \).\label{lt.10}
 \end{equation}
for all $n$, $x_{1},\ldots x_{n}\in S$ and $\la_{1}, \ldots, \la_{n}$ sufficiently small, since we can always take $\la_{1}=\la_{2}=0$. By the independence of $X$ and $G$ this is equivalent to showing that
\begin{equation}
Q^{x_{1},x_{2}}\(e^{\sum_{j=1}^{n}\la_{j}\,L^{ x_{j}}_{\ff}}  \)=\frac{E_{G} \(G_{ x_{1}}G_{ x_{2}}\, e^{\sum_{j=1}^{n}\la_{j}{1 \over 2}G^{2}_{ x_{j}}} \)}{E_{G} \( e^{\sum_{j=1}^{n}\la_{j}{1 \over 2}G^{2}_{ x_{j}}} \)}.\label{lt.11}
\end{equation}
Differentiating (\ref{lt.4}) with respect to $u_{1}, u_{2}$ and then setting all $u_{j}=0$ for the numerator and using (\ref{lt.4}) with  all $u_{j}=0$ for the denominator we see that
\begin{equation}
\frac{E_{G} \(G_{ x_{1}}G_{ x_{2}}\, e^{\sum_{j=1}^{n}\la_{j}{1 \over 2}G^{2}_{ x_{j}}} \)}{E_{G} \( e^{\sum_{j=1}^{n}\la_{j}{1 \over 2}G^{2}_{ x_{j}}} \)}=\bar C_{1,2}.\label{lt.12}
\end{equation}

To evaluate the left hand side of (\ref{lt.11}) it is useful to introduce the atomic measure on $S$
\begin{equation}
\nu=\sum_{j=1}^{n}\la_{j}\de_{x_{j}}\label{lt.13}
\end{equation}
and to write $\sum_{j=1}^{n}\la_{j}\,L^{ x_{j}}_{\ff} =\int L^{ x}_{\ff}\,d\nu(x)$. With this notation we obtain from (\ref{mp.14})
  \begin{eqnarray}
 &&Q^{x_{1},x_{2}}\( \(\int L^{ x}_{\ff}\,d\nu(x)\)^{k} \) 
 \label{lt.15}\\
 &&  = k!\int  u(x_{1},y_{1}) u(y_{1},y_{2})\cdots   u(y_{k-1},y_{k})u(y_{k},x_{2})\prod_{j=1}^{k}\,d\nu (y_{j}) \nonumber\\
 &&  = k! \,\,(\,\,\overset{k}{\overbrace{C\La \cdots C\La} }\,\, C\,\,)_{1,2}=k! ((C\La)^{k}C)_{1,2}.\nonumber
 \end{eqnarray}
 Thus 
 \begin{equation}
 Q^{x_{1},x_{2}}\(e^{\sum_{j=1}^{n}\la_{j}\,L^{ x_{j}}_{\ff}}  \)=\sum_{k=0}^{\ff}((C\La)^{k}C)_{1,2}=\bar C_{1,2}.\label{lt.16}
 \end{equation}
 \qed

\subsection{Another proof of the Eisenbaum Isomorphism theorem}

As in the last subsection, it suffices to prove that
\begin{equation}
P^{x_{1}}\(e^{\sum_{j=1}^{n}\la_{j}\,L^{ x_{j}}_{\ff}}  \)=\frac{E_{G} \(\(1+{G_{ x_{1}} \over s}\) \, e^{\sum_{j=1}^{n}\la_{j}{1 \over 2}(G_{ x_{j}}+s)^{2}} \)}{E_{G} \( e^{\sum_{j=1}^{n}\la_{j}{1 \over 2}(G_{ x_{j}}+s)^{2}} \)}.\label{lt.17}
\end{equation}
for all $n$, $x_{1},\ldots x_{n}\in S$ and $\la_{1}, \ldots, \la_{n}$ sufficiently small. We can write the right hand side as
\begin{equation}
1+\frac{E_{G} \( G_{ x_{1}}  \, e^{\sum_{j=1}^{n}s\la_{j}G_{ x_{j}} +\la_{j}{1 \over 2}G_{ x_{j}}^{2}} \)}{sE_{G} \( e^{\sum_{j=1}^{n}s\la_{j}G_{ x_{j}} +\la_{j}{1 \over 2}G_{ x_{j}}^{2}} \)}=1+ \sum_{j=1}^{n}\bar C_{1,j}\la_{j},\label{lt.18}
\end{equation}
where the last equality comes from differentiating (\ref{lt.4}) with respect to $u_{1}$ and then setting all $u_{j}=s$ for the numerator, and setting all $u_{j}=s$ for the denominator. We can then rewrite
\begin{equation}
1+ \sum_{j=1}^{n}\bar C_{1,j}\la_{j}=1+\sum_{j=1}^{n}\lc(I-C\La )^{-1}C\La\rc_{1,j}=
1+  \sum_{j=1}^{n}(\sum_{k=1}^{\ff}(C\La)^{k})_{1,j}. \label{lt.18a}
\end{equation}

Using again the notation (\ref{lt.13}), we obtain from (\ref{mp.12})
  \begin{eqnarray}
 &&P^{x_{1}}\( \(\int L^{ x}_{\ff}\,d\nu(x)\)^{k} \) 
 \label{lt.19}\\
 &&  = k!\int  u(x_{1},y_{1}) u(y_{1},y_{2})\cdots   u(y_{k-1},y_{k})\prod_{j=1}^{k}\,d\nu (y_{j}) \nonumber\\
 &&  = k! \,\,\sum_{j=1}^{n}(\,\,\overset{k}{\overbrace{C\La \cdots C\La} })_{1,j}=k! \sum_{j=1}^{n}((C\La)^{k})_{1,j}.\nonumber
 \end{eqnarray}
Hence
\begin{equation}
P^{x_{1}} \(e^{\sum_{j=1}^{n}\la_{j}\,L^{ x_{j}}_{\ff}}  \)=\sum_{j=1}^{n}(\sum_{k=0}^{\ff}(C\La)^{k})_{1,j}\label{lt.20}
\end{equation}
which is the same as (\ref{lt.18a}), since $\sum_{j=1}^{n}((C\La)^{0})_{1,j}=\sum_{j=1}^{n}I_{1,j}=1.$\qed

\subsection{Another proof of the generalized second Ray-Knight theorem}

As in the last two subsections, it suffices to prove that for all $t$
\begin{equation}
P^{0}\(e^{\sum_{j=1}^{n}\la_{j}\,L^{ x_{j}}_{\tau(t)}}  \)=\frac{E_{\eta} \( e^{\sum_{j=1}^{n}\la_{j}{1 \over 2}(\eta_{ x_{j}}+\sqrt{2t})^{2}} \)}{E_{\eta} \( e^{\sum_{j=1}^{n}\la_{j}{1 \over 2}\eta_{ x_{j}}^{2}} \)}.\label{lt.21}
\end{equation}
for all $n$, $x_{1},\ldots x_{n}\in S$ and $\la_{1}, \ldots, \la_{n}$ sufficiently small. Let $C_{0}$ denote the covariance matrix of $(\eta_{1}, \ldots, \eta_{n})$ and let $\bar 1=(1,\ldots,1)$, that is, the vector in $R^{n}$ with all componets equal to $1$. 
Using (\ref{lt.8}) with  all $u_{j}=\sqrt{2t}$ for the numerator, and   all $u_{j}=0$ for the denominator
we obtain
\begin{equation}
\frac{E_{\eta} \( e^{\sum_{j=1}^{n}\la_{j}{1 \over 2}(\eta_{ x_{j}}+\sqrt{2t})^{2}} \)}{E_{\eta} \( e^{\sum_{j=1}^{n}\la_{j}{1 \over 2}\eta_{ x_{j}}^{2}} \)}=e^{t(\bar 1,\La \bar 1)}
e^{t(\bar 1,\La \bar C_{0}\La \bar 1)}=e^{t\(\bar 1,\La\(\sum^{\ff}_{k=0} (C_{0}\La)^{k} \bar 1\)\)},\label{lt.22}
\end{equation}
where the last equality used (\ref{lt.9}).

Using once again the notation (\ref{lt.13}), we obtain from (\ref{rk.5})
  \begin{eqnarray}
 &&P^{0}\( \(\int L^{ x}_{\tau(t)}\,d\nu(x)\)^{k} \) 
 \label{lt.23}\\
 &&  =\sum_{m=1}^{k}\sum_{\stackrel{\mbox{\scriptsize unordered}}{B_{1}\cup\cdots\cup B_{m}=[1,k]}} \,\,t^{m} \int \prod_{j=1}^{m} \mbox{ ch}_{0} (B_{j})\prod_{l=1}^{k}\,d\nu (y_{l}). \nonumber
 \end{eqnarray}
 Hence if we set
 \begin{equation}
 h(k)=\int u_{T_{0}}(y_{1},  y_{2})\cdots u_{T_{0}}(y_{k-1},  y_{k})\prod_{l=1}^{k}\,d\nu (y_{l})=
\(\bar 1,\La  (C_{0}\La)^{k-1} \bar 1\)\label{lt.25}
 \end{equation}
we see that  for any partition $B_{1}\cup\cdots\cup B_{m}=[1,k]$
\begin{equation}
\int \prod_{j=1}^{m} \mbox{ ch}_{0} (B_{j})\prod_{l=1}^{k}\,d\nu (y_{l}) =\prod_{j=1}^{m}|B_{j}|!\,h(|B_{j}|). \label{lt.26}
\end{equation}
 Since there are ${1 \over m!}{k \choose k_{1}\cdots k_{m}}$ ways to partition $k$ objects into $m$ unordered  subsets of size $k_{1},\cdots, k_{m}$ we see that
 \bea
P^{0}\( \(\int L^{ x}_{\tau(t)}\,d\nu(x)\)^{k} \) &=& \sum_{m=1}^{k}{t^{m} \over m!}\sum_{k_{1}+\cdots+ k_{m}=k}{k \choose k_{1}\cdots k_{m}}\prod_{j=1}^{m}k_{j}!h(k_{j} )\nn\\
&=&   k!\sum_{m=1}^{k}{t^{m} \over m!} \sum_{k_{1}+\cdots+ k_{m}=k} \prod_{j=1}^{m} h(k_{j} )       \label{lt.27}
 \eea
 Hence
 \be
P^{0}\( e^{\int L^{ x}_{\tau(t)}\,d\nu(x)}   \) =\sum_{m=0}^{\ff}{t^{m} \over m!}\(\sum_{j=1}^{\ff}h(j) \)^{m}=e^{t \sum_{j=1}^{\ff}h(j)}.
 \label{lt.28}
 \ee
In view of (\ref{lt.25}), this gives (\ref{lt.22}).

\subsection{Yet another proof of the generalized second Ray-Knight theorem using excursion theory}\label{subsec-yet}

It follows from (\ref{et.1}) and the master formula (\ref{pp.3}) that for $\de$ small
\begin{equation}
P^0 \(e^{\de \int L^{ x}_{\tau(t)}\,d\nu(x)}\)=\exp \(t{\bf  n}\(e^{\de \int L^{ x}_{\ff}\,d\nu(x)}-1\)\),\label{4.4}
\end{equation}
and it follows from (\ref{et.111}) and (\ref{lt.26}) that 
\be{\bf  n}\(e^{\de \int L^{ x}_{\ff}\,d\nu(x)}-1\)=\sum_{n=1}^{\ff}\de^{n } \,h( n).
\label{4.5}
\ee
This completes the proof of (\ref{lt.28}) which we have seen is sufficient to prove our theorem.

\subsection{Another proof of the interlacement Isomorphism theorem}\label{sec-intsymlap}

Because everything is additive in $\nu$ we can write (\ref{q2.0}) as 
\be
{P}_{\mathcal{I}_{t}}\times P_{G}\(e^{\de \wt L^\nu_\ff+{\de\over 2}:G^{2}:(\nu)}\)= P_{G}\(e^{{ \de\over 2}:G^{2}:(\nu)+\de\sqrt{2 t}G_{\nu}+\de t  |\nu|}\)\label{q4.1}
\ee
for $\de$ small. Equivalently, we show that
\be
{P}_{\mathcal{I}_{t}}\(e^{\de  \wt L^\nu_\ff}\)= { P_{G}\(e^{{ \de\over 2}:G^{2}:(\nu)+\de\sqrt{2 t}G_{\nu}+\de t  |\nu|}\)\over P_{G}\(e^{{ \de\over 2}: G^2:(\nu)}\) }.\label{q4.2}
\ee

 (\ref{q.4}) shows that
\begin{equation}
\mbox{$ \bf{P}_{m}$}\(\(L_{\ff}^{\nu}\)^{k}\)=k!\int  \prod_{j=1}^{k-1} u(y_{j},y_{j+1})\prod_{j=1}^{k}\nu(dy_{j}).\label{q.8}
\end{equation}
Given (\ref{q.8}) and our  use of the master formula in subsection \ref{subsec-yet}
it suffices to show that 
\bea
&&
 { P_{G}\(e^{{ \de\over 2}:G^{2}:(\nu)+\de\sqrt{2 t}G_{\nu}+\de t  |\nu|}\)\over P_{G}\(e^{{ \de\over 2}: G^2:(\nu)}\) }\label{q4.3}\\
  &&\hspace{1 in}=\exp \( t\(\sum_{n=1}^{\ff}\de^{n }\int \prod_{j=1}^{n-1} u(x_{j},x_{j+1})\prod_{j=1}^{n}\nu(dx_{j})\)\).\nn
\eea
To see this, we first note that using (\ref{a.3a}), the Gaussian moment formula and the monotone convergence theorem  we have 
\begin{equation}
P_{G}\(e^{{ \de\over 2}:\(G\)^2:(\nu)}\)=\lim_{\ep\to 0}P_{G_{\de}}\(e^{{ \de\over 2}\int \( G^{2}_{x,\ep}-E\(G^{2}_{x,\ep}\)   \) \,d\nu(x)}\)\label{q4.4}
\end{equation}
and
\bea
&&
P_{G}\(e^{{ \de\over 2}:G^{2}:(\nu)+\de\sqrt{2 t}G_{\nu}+\de t|\nu|}\)\label{q4.5}\\
&&=\lim_{\ep\to 0}P_{G_{\de}}\(e^{{ \de\over 2}\int \((G_{x,\ep}+\sqrt{2 t})^{2} -E\(G^{2}_{x,\ep}\)   \) \,d\nu(x)}\)\nn
\eea
Therefore
\begin{eqnarray}
&& {P_{G}\(e^{{ \de\over 2}:G^{2}:(\nu)+\de\sqrt{2 t}G_{\nu}+\de t  |\nu|}\)\over P_{G}\(e^{{ \de\over 2}: G^2:(\nu)}\) }
\label{q4.6}\\
&&=\lim_{\ep\to 0}   {P_{G_{\de}}\(e^{{ \de\over 2}\int \((G_{x,\ep}+\sqrt{2 t})^{2} -E\(G^{2}_{x,\ep}\)   \) \,d\nu(x)}\) \over P_{G_{\de}}\(e^{{ \de\over 2}\int \( G^{2}_{x,\ep}-E\(G^{2}_{x,\ep}\)   \) \,d\nu(x)}\)}\nonumber\\
&&=\lim_{\ep\to 0}   {P_{G_{\de}}\(e^{{ \de\over 2}\int  (G_{x,\ep}+\sqrt{2 t})^{2}      \,d\nu(x)}\) \over P_{G_{\de}}\(e^{{ \de\over 2}\int   G^{2}_{x,\ep}     \,d\nu(x)}\)}\nonumber\\
&&=\lim_{\ep\to 0}  \exp \( t \(\sum_{n=1}^{\ff}\de^{n }\int \prod_{j=1}^{n-1} u_{\ep}(x_{j},x_{j+1})\prod_{j=1}^{n}\nu(dx_{j})\)\) \nonumber
\end{eqnarray}
as in (\ref{lt.22}).
Using the monotone convergence theorem this is
\begin{equation}
 =\exp \( t \(\sum_{n=1}^{\ff}\de^{n }\int \prod_{j=1}^{n-1} u(x_{j},x_{j+1})\prod_{j=1}^{n}\nu(dx_{j})\)\).\label{}
\end{equation}
This completes the proof of (\ref{q4.3}) and hence of (\ref{q2.0}).\qed

\def\noopsort#1{} \def\printfirst#1#2{#1}
\def\singleletter#1{#1}
            \def\switchargs#1#2{#2#1}
\def\bibsameauth{\leavevmode\vrule height .1ex
            depth 0pt width 2.3em\relax\,}
\makeatletter
\renewcommand{\@biblabel}[1]{\hfill#1.}\makeatother
\newcommand{\bysame}{\leavevmode\hbox to3em{\hrulefill}\,}

 \def\wh{\widehat}
\def\ol{\overline}


\begin{thebibliography}{10}

 

\bibitem{BES}
Bass, R.~F., Eisenbaum, N.,  and  Shi, Z. 
\newblock The most visited sites of symmetric stable processes.
\newblock {\em Prob.\ Theory\ Related\ Fields}, {\em 116}, (2000), 391--404.


\bibitem{Bertoin} J.  Bertoin {\em Random fragmentation and coagulation processes}, Cambridge University Press, New
York,  (2006).  




 \bibitem{B}
R. M. Blumenthal, (1992)
\newblock{\em Excursions of Markov Processes},
\newblock Birkhauser, Boston.




\bibitem {DM3}
C. Dellacherie,  and P.-A. Meyer,  (1987).
\newblock {\em Probabilities et Potential, Chapitres XII a XVI}.
\newblock Paris: Hermann.

\bibitem {DM4}
C. Dellacherie,  B. Maisonneuve and P.-A. Meyer,  (1992).
\newblock {\em Probabilities et Potential, Chapitres XVII a XXIV}.
\newblock Paris: Hermann.


\bibitem{DLP}
J. Ding, J. Lee and Y. Peres.
Cover times, blanket times, and majorizing measures.  
Annals of Math 175(3) : 1409-1471 (2012), conference version at STOC (2011).

\bibitem{D1}
J. Ding. 
Asymptotics of cover times via Gaussian free fields: bounded-degree graphs and general trees. Annals of Probability, to appear. 

\bibitem{D2}
J. Ding.
On cover times for 2D lattices, {\it  EJP},\, to appear.

\bibitem{DZ1}
J. Ding and O. Zeitouni.
A sharp estimate for cover times on binary trees, {\it SPA},\, to appear.

\bibitem{DZ2}
J. Ding and O. Zeitouni.
Extreme values for two-dimensional discrete Gaussian free field. Annals of Probability, to appear.


\bibitem{Dynkin83}
Dynkin, E.~B.  
\newblock Local times and quantum fields.
\newblock In {\em Seminar on Stochastic Processes}, volume~7 of {\em Progress
  in Probability}, (1983),  (pp.\ 64--84). Boston: Birkh\"auser.

\bibitem{Dynkin84}
Dynkin, E.~B. 
\newblock {G}aussian and non-{G}aussian random fields associated with {M}arkov
  processes.
\newblock {\em J.~Fcnl.\ Anal.}, {\em 55}, (1984), 344--376.




 \bibitem {five}
N. Eisenbaum, H. Kaspi, M. Marcus, J. Rosen and Zhan Shi, 
A Ray-Knight theorem for symmetric Markov processes,
 {\it Ann. Probab.},\,
 {\bf 28}\, (2000), 1781-1796.

 
 
 


\bibitem {EK}
N. Eisenbaum  and  H. Kaspi,
\newblock On permanental processes, {\em Stochastic Processes and their Applications},
{\em 119},  (2009),  1401-1415.


\bibitem {FPY}
P. Fitzsimmons, J. Pitman,  and M. Yor,  
\newblock Markovian bridges: construction, Palm interpretation, and splicing.  
{\em Seminar on Stochastic Processes}, 1992, E. Cinlar and K.L. Chung and M.J. Sharpe editors, 101-134, BirkhŠuser, Boston (1993).

\bibitem {FR}
P. Fitzsimmons and J.~Rosen,  
\newblock Markovian loop soups: permanental processes and isomorphism theorems.  \, {\it Electron. J. Probab.},\, Volume 19 (2014), no. 60, 1-30.
http://arxiv.org/pdf/1211.5163.pdf 

\bibitem{K} J. F. C. Kingman, {\em Poisson Processes}, Oxford Studies in Probability, Clarendon Press, Oxford, (2002).  

\bibitem{LL} G. Lawler and V. Limic, {\em Random Walk: A Modern Introduction}, Cambridge University Press, New
York,  (2009).

\bibitem{LF}
G. Lawler and J. Trujillo Ferreis, \newblock Random walk loop soup, \newblock {\em TAMS } 359  (2007), 565--588.

\bibitem{LW}
G. Lawler and W. Werner,  \newblock The Brownian loop soup, \newblock {\em PTRF} 44 (2004), 197--217.



 
 \bibitem{Le Jan} Y. Le Jan, \newblock
  Markov loops and renormalization, {\em Ann. Probab.},
{\bf 38}  (2010),  1280--1319. 

\bibitem{Le Jan1} Y. Le Jan, \newblock
{\em Markov paths, loops and fields. }   \'{E}cole d'\'{E}t\'{e} de Probabilit\'{e}s de Saint-Flour XXXVIII - 2008. Lecture Notes in Mathematics 2026. 
Springer-Verlag, Berlin-Heidelberg, (2011).

\bibitem{LMR} Y. Le Jan, M. B.  Marcus and J.~Rosen,\newblock
{\em Permanental fields,  loop soups and continuous additive functionals.},\, {\it Ann. Probab.},\, to appear. http://arxiv.org/pdf/1209.1804.pdf 

\bibitem{LMR2} Y. Le Jan, M. B.  Marcus and J.~Rosen,
\newblock  Intersection local times, loop soups and permanental Wick powers.
  http://arxiv.org/pdf/1308.2701.pdf



\bibitem{book} M. B.  Marcus and J.~Rosen, {\em Markov Processes,
Gaussian Processes and Local Times}, Cambridge University Press, New
York,  (2006).

 
  
\bibitem{MRperm} M. B. Marcus and J.~Rosen, A sufficient condition for the continuity of permanental processes with
 applications to local times of Markov processes,\, {\it Ann. Probab.}, \, 41, (2013), 671--698.\\ http://arxiv.org/pdf/1005.5692.pdf 
 
 \bibitem{MRcontperm} M. B. Marcus and J.~Rosen,   Continuity conditions for a class of second order permanental chaoses,\, High Dimensional Probability VI: the Banff volume, Progress in Probability, 66 (2013), 229-245, Springer, Basel.
 

\bibitem{RW} 
  L. C. G. Rogers and D. Williams,
  {\em Diffusions, {M}arkov Processes, and Martingales. Volume Two:
 Ito Calculus},
  Cambridge University Press,
 Cambridge,  (2000).
 

\bibitem{ILTInt} J.~Rosen, Intersection local times for interlacements.  {\it Stochastic Processes and their Applications},\,  Volume 124, Issue 5, May 2014, Pages 1849-1880.   arxiv.org/pdf/1308.3469.pdf 

  
   \bibitem{S} M. Sharpe, {\em General theory of Markov
processes},  Acad. Press, New York, (1988).


\bibitem{Sz1} A.-S. Sznitman,
\newblock  Topics in occupation times and Gaussian free fields.
 {\it Zurich Lectures in Advanced Mathematics}, EMS, Zurich, 2012.
 
 \bibitem{Sz2} A.-S. Sznitman,
\newblock  An isomorphism theorem for random interlacements.


\bibitem{VJ}  D.  Vere-Jones,
     {\em Alpha-permanents},
    New Zealand J. of Math.,  (1997), 26, 125--149.
    
    
     




 

 
 
 

 



\end{thebibliography}
\end{document}